\documentclass[12pt,a4paper]{article}

\usepackage{amsfonts,amsmath,amsthm,amssymb,amsopn}
\usepackage{graphicx}
\usepackage{epstopdf}
\usepackage[export]{adjustbox}
\usepackage{algorithmic}
\usepackage{enumitem}

\usepackage{bm,color}
\usepackage{xcolor}
\usepackage{stmaryrd}
\usepackage[T1]{fontenc}
\usepackage{float, caption, subcaption}
\usepackage{booktabs}
\usepackage{multirow}
\usepackage{siunitx}
\usepackage{pifont}
\usepackage{hyperref}
\usepackage{cleveref}
\usepackage{bbold}

\newcommand\bF{\bm{F}}

\newcommand\bv{\bm{v}}

\newcommand\bU{\bm{U}}
\newcommand\bx{\bm{x}}

\newcommand\bB{\bm{B}}

\newcommand\dd{\mathrm{d}}
\newcommand\pd[2]{\frac{\partial {#1}}{\partial {#2}}}

\newcommand\abs[1]{\left\lvert #1 \right\rvert}
\newcommand\norm[1]{\left\lVert #1 \right\rVert}

\newcommand\xr{i+\frac12}
\newcommand\xl{i-\frac12}
\newcommand\yr{j+\frac12}
\newcommand\yl{j-\frac12}

\newcommand\divB{\nabla\cdot{\bB}}
\newcommand\vB{\bv\cdot{\bB}}

\theoremstyle{plain}% default
\newtheorem{definition}{Definition}[section]

\theoremstyle{definition}
\newtheorem{lemma}{Lemma}[section]

\newtheorem{remark}{Remark}[section]
\newtheorem{example}{Example}[section]

\crefname{equation}{Equation}{Equations}
\crefname{figure}{Figure}{Figures}
\crefname{table}{Table}{Tables}
\crefname{example}{Example}{Examples}
\crefname{section}{Section}{Sections}

\renewcommand{\title}{Active flux for ideal magnetohydrodynamics: A positivity-preserving scheme with the Godunov-Powell source term}

\newcommand{\authorOne}{Junming Duan\footnote{Institute of Mathematics, University of W\"urzburg, Emil-Fischer-Stra\ss e 40, 97074 W\"urzburg, Germany, junming.duan@uni-wuerzburg.de}}
\newcommand{\authorTwo}{Praveen Chandrashekar\footnote{Centre for Applicable Mathematics, TIFR, 560065, Bangalore, India, praveen@tifrbng.res.in}}
\newcommand{\authorThree}{Christian Klingenberg\footnote{Institute of Mathematics, University of W\"urzburg, Emil-Fischer-Stra\ss e 40, 97074 W\"urzburg, Germany, christian.klingenberg@uni-wuerzburg.de}}
\begin{document}

\begin{center} \Large
\title

\vspace{1cm}

\date{}
\normalsize

\authorOne, \authorTwo, \authorThree
\end{center}

\begin{abstract}

% #############################################################################################################################
The Active Flux (AF) is a compact, high-order finite volume scheme that allows more flexibility by introducing additional point value degrees of freedom at cell interfaces.
This paper proposes a positivity-preserving (PP) AF scheme for solving the ideal magnetohydrodynamics, where the Godunov-Powell source term is employed to deal with the divergence-free constraint.
For the evolution of the cell average, apart from the standard conservative finite volume method for the flux derivative, the nonconservative source term is built on the quadratic reconstruction in each cell, which maintains the compact stencil in the AF scheme.
For the point value update, the local Lax-Friedrichs (LLF) flux vector splitting is adopted for the flux derivative, originally proposed in [Duan, Barsukow, and Klingenberg, SIAM Journal on Scientific Computing, 47(2), A811--A837, 2025], and a central difference is used to discretize the divergence in the source term.
A parametrized flux limiter and a scaling limiter are presented to preserve the density and pressure positivity by blending the AF scheme with the first-order PP LLF scheme with the source term.
To suppress oscillations, a new shock sensor considering the divergence error is proposed, which is used to compute the blending coefficients for the cell average.
Several numerical tests are conducted to verify the third-order accuracy, PP property, and shock-capturing ability of the scheme.
The key role of the Godunov-Powell source term and its suitable discretization in controlling divergence error is also validated.

% #############################################################################################################################



\end{abstract}

% #############################################################################################################################
\section{Introduction}\label{sec:introduction}
This paper is concerned with solving the two-dimensional ideal magnetohydrodynamics (MHD), which in conservative form reads
\begin{equation}\label{eq:2d_mhd}
	\bU_t + \bF_1(\bU)_x + \bF_2(\bU)_y = 0,
\end{equation}
where $\bU=(\rho, \rho\bv^\top, \bB^\top, E)^\top$ is the vector of conservative variables,
with the velocity vector $\bv=(v_1, v_2, v_3)^\top$ and magnetic field vector $\bB=(B_1, B_2, B_3)^\top$,
the total energy $E = \frac12\rho\norm{\bv}^2 + \rho e + \frac12\norm{\bB}^2$,
and $\rho e$ is the internal energy.
The fluxes in the $x$- and $y$-directions are defined as
\begin{equation*}
	\bF_1=
	\begin{pmatrix}
		\rho v_1\\ \rho v_1^2 - B_1^2 + p_t \\ \rho v_1v_2 - B_1B_2\\ \rho v_1v_3 - B_1B_3\\ 0\\ v_1B_2-B_1v_2\\ v_1B_3-B_1v_3\\ (E+p_t)v_1 - B_1(\vB) \\
	\end{pmatrix},\quad
	\bF_2=
	\begin{pmatrix}
		\rho v_2\\ \rho v_1v_2 - B_1B_2 \\ \rho v_2^2 - B_2^2 + p_t \\ \rho v_2v_3 - B_2B_3\\ v_2B_1-B_2v_1\\ 0\\ v_2B_3-B_2v_3\\ (E+p_t)v_2 - B_2(\vB) \\
	\end{pmatrix}.
\end{equation*}
Here the total pressure $p_t = p + p_m$ consists of the fluid pressure $p$ and magnetic pressure $p_m = \frac12\norm{\bB}^2$.
To close the system \eqref{eq:2d_mhd}, this paper considers the equation of state (EOS) for the perfect gas $p = (\gamma-1)\rho e$,
with the adiabatic index $\gamma$.
The physical solutions should satisfy the divergence-free constraint on the magnetic field
\begin{equation*}
	\divB = \pd{B_1}{x} + \pd{B_2}{y} = 0.
\end{equation*}

For the numerical solutions of the MHD equations \eqref{eq:2d_mhd}, one needs to carefully deal with the divergence-free constraint, otherwise, large divergence errors may lead to nonphysical features or numerical instabilities \cite{Brackbill_1980_effect_JCP,Balsara_1999_staggered_JCP,Toth_2000_$B0$_JCP}.
Many works have focused on this issue, e.g., the projection method \cite{Brackbill_1980_effect_JCP},
the constrained transport (CT) method and its variants \cite{Evans_1988_Simulation_AJ,Gardiner_2005_unsplit_JCP,Londrillo_2004_divergence_JCP},
the eight-wave formulation of the MHD equations \cite{Powell_1994_approximate, Powell_1999_solution_JCP} based on the Godunov-Powell source term \cite{Godunov_1972_Symmetric,Powell_1994_approximate}, the hyperbolic divergence cleaning method \cite{Dedner_2002_Hyperbolic_JCP},
the locally divergence-free discontinuous Galerkin (DG) method \cite{Li_2005_Locally_JSC},
the globally divergence-free central DG method \cite{Li_2011_Central_JCP}, etc.
By adding the Godunov-Powell source term \cite{Godunov_1972_Symmetric,Powell_1994_approximate} to the conservative MHD equations \eqref{eq:2d_mhd},
the modified system reads
\begin{equation}\label{eq:2d_mhd_gp}
	\bU_t + \bF_1(\bU)_x + \bF_2(\bU)_y = - (\divB)\bm{\Psi},
\end{equation}
where
\begin{equation}\label{eq:gp_src}
	\bm{\Psi} = (0, \bB, \bv, \vB)^\top.
\end{equation}
Such a source term makes the system nonconservative but introduces many advantages, e.g., the modified system is Galilean invariant, and can be symmetrized by the entropy pair \cite{Godunov_1972_Symmetric}, which leads to the entropy stable schemes for the MHD \cite{Chandrashekar_2016_Entropy_SJNA,Winters_2016_Affordable_JCP,Liu_2018_Entropy_JCP}.
One can also verify that the divergence satisfies the following transport equation,
\begin{equation*}
	\dfrac{\partial}{\partial t}\left(\dfrac{\divB}{\rho}\right)
	+ \bv\cdot\nabla \left(\dfrac{\divB}{\rho}\right) = 0,
\end{equation*}
which means divergence error may be advected away by the flow \cite{Powell_1994_approximate}, instead of accumulating and causing instabilities.

The design of so-called positivity-preserving (PP) numerical methods that maintain the positivity of density and pressure is also very important for numerical stability.
To address this issue, several techniques have been proposed \cite{Balsara_1999_Maintaining_JCP,Janhunen_2000_positive_JCP,Balsara_2012_Self_JCP}.
The PP Riemann solver based on relaxation was constructed in \cite{Bouchut_2007_Multiwave_NM,Bouchut_2010_Multiwave_NM},
and PP schemes based on that were studied in \cite{Waagan_2009_positive_JCP, Waagan_2011_robust_JCP}.\
It was shown in \cite{Bouchut_2010_Multiwave_NM} that the Godunov-Powell source term is important in the design of PP schemes for the multi-dimensional MHD.
The PP DG and central DG schemes based on the scaling limiter \cite{Zhang_2010_positivity_JCP} were proposed in \cite{Cheng_2013_Positivity_JCP}.
The PP finite difference methods were developed in \cite{Christlieb_2015_Positivity_SJSC,Christlieb_2016_High_JCP} by using the parametrized flux limiter \cite{Xu_2014_Parametrized_MC}.
Taking inspiration from transforming nonlinear constraints of the admissible state set into linear ones by adding auxiliary variables \cite{Wu_2017_Admissible_MMMAS},
the first-order Lax-Friedrichs (LF) scheme with suitable viscosity and a discrete Godunov-Powell source term was rigorously proved to be PP by Wu in \cite{Wu_2018_Positivity_SJNA},
and based on that, the provably high-order PP DG schemes were proposed in \cite{Wu_2018_provably_SJSC} with a suitable high-order discretization of the Godunov-Powell source term.
It is also demonstrated in \cite{Wu_2018_provably_SJSC,Wu_2019_Provably_NM} that the Godunov-Powell source term helps to eliminate the effect of divergence error on the PP property.
Subsequently, more works on the design of high-order PP schemes for the MHD equations were presented, including but not restricted to \cite{Wu_2019_Provably_NM,Ding_2024_New_SJSC,Liu_2025_Structure_JCP}.

The active flux (AF) method is a compact finite volume method \cite{Eymann_2011_Active_InCollection, Eymann_2011_Active_InProceedings,Eymann_2013_Multidimensional_InCollection,Roe_2017_Is_JSC},
with inspiration from \cite{VanLeer_1977_Towards_JCP}.
It simultaneously evolves cell averages and additional degrees of freedom (DoFs), chosen as point values at cell interfaces like the continuous finite element method.
Thanks to this continuity of the point values across the cell interface, the AF method does not need Riemann solvers (unlike Godunov methods) for the evolution of the cell average.
The AF methods can be roughly divided into two classes based on the evolution of the point value.
The original ones evolve the cell average through Simpson's rule for flux quadrature in time,
and employ exact or approximate evolution operators to evolve the point values and to obtain the numerical solutions at the flux quadrature points.
Examples are exact evolution operators for linear equations \cite{Barsukow_2019_Active_JSC,Fan_2015_Investigations_InCollection,Eymann_2013_Multidimensional_InCollection, VanLeer_1977_Towards_JCP},
$p$-system \cite{Fan_2017_Acoustic},
and approximate evolution operators for Burgers' equation \cite{Eymann_2011_Active_InCollection,Eymann_2011_Active_InProceedings,Roe_2017_Is_JSC,Barsukow_2021_active_JSC},
the compressible Euler equations in one spatial dimension \cite{Eymann_2011_Active_InCollection,Helzel_2019_New_JSC,Barsukow_2021_active_JSC},
multidimensional Euler equations \cite{Fan_2017_Acoustic},
and hyperbolic balance laws \cite{Barsukow_2021_Active_SJSC, Barsukow_2023_Well_CAMC}, etc.
The method of bicharacteristics was used for the derivation of truly multidimensional approximative evolution operators \cite{Chudzik_2024_Active_JSC}.
The other so-called generalized, or semi-discrete AF methods adopt a method of lines,
where the evolution of the cell average and point value is written in semi-discrete form and integrated in time by using Runge-Kutta methods.
Examples of this approach are \cite{Abgrall_2023_Combination_CAMC, Abgrall_2023_Extensions_EMMNA, Abgrall_2025_semi_JSC, Abgrall_2024_Bound} based on Jacobian splitting (JS) and \cite{Duan_2025_Active_SJSC} based on flux vector splitting (FVS).
The AF method is superior to standard finite volume methods due to its structure-preserving property.
For example, it preserves the vorticity and stationary states for multi-dimensional acoustic equations \cite{Barsukow_2019_Active_JSC},
and it is naturally well-balanced for acoustics with gravity \cite{Barsukow_2021_Active_SJSC}.

This paper proposes a PP AF scheme for solving the ideal MHD equations, where the Godunov-Powell source term is employed to deal with the divergence-free constraint.
For the discretization of the flux derivative in the point value update, we use the local Lax-Friedrichs (LLF) FVS following the previous work \cite{Duan_2025_Active_SJSC}, which shows better performance for strong discontinuities.
Our main novelty and contributions in this paper are as follows.
\begin{itemize}
	\item We construct suitable discretizations for the nonconservative source term to achieve a stable AF scheme.
	For the evolution of cell average, the source term is discretized using the $3\times 3$ Gauss-Lobatto quadrature rule, where the discrete divergence is easily computed based on the quadratic reconstruction in each cell, thus such a discretization only depends on the DoFs in the current cell, and maintains the compactness of the AF scheme.
	For the point value update, a central difference is used to discretize the divergence in the source term, built on the same spatial stencil as the original AF scheme \cite{Abgrall_2025_semi_JSC}.
	Numerical tests will show that the inclusion of the Godunov-Powell source term and our discretization can control the divergence error.
	\item To design the PP AF scheme, we borrow the idea of blending the high-order AF scheme and the first-order PP LLF scheme from \cite{Duan_2025_Active_SJSC}, and take advantage of the PP property of the first-order LLF scheme, which was rigorously proved in \cite{Wu_2019_Provably_NM}.
	Different from \cite{Duan_2025_Active_SJSC}, a parametrized flux limiter \cite{Xu_2014_Parametrized_MC,Christlieb_2015_Positivity_SJSC} is adopted for the cell average as the intermediate state defined in \cite{Duan_2025_Active_SJSC} may not be PP for the MHD.
	Our PP limiting for the cell average consists of two steps: the source term is blended first and then the numerical flux.
	A scaling limiter, as the one in \cite{Duan_2025_Active_SJSC}, is also presented to preserve the PP property for the point value update.
	Thus, our AF scheme is PP for both the cell average and point value.
	\item To suppress oscillations, a new shock sensor is proposed to be used in the blending for the cell average.
	We take into account the magnetic pressure and also divergence error,	where the latter indicates the nonsmooth regions in the magnetic field.
	Additionally, we also limit the discretization for the source term based on the blending coefficients at cell edges.
	Several numerical examples, including the rotor problem, blast problem, and high Mach number jets in a strongly magnetized medium, will be used to demonstrate the ability of the shock sensor-based limiting.
\end{itemize}

The remainder of this paper is structured as follows.
\Cref{sec:2d_af_scheme} constructs the 2D AF scheme based on the LLF FVS for the point value update, and suitable discretizations for the Godunov-Powell source term.
\Cref{sec:2d_limitings} presents the 2D PP limitings and also a limiting for suppressing oscillations based on a new shock sensor.
Numerical tests are conducted in \Cref{sec:results} to experimentally demonstrate the accuracy, PP property, and shock-capturing ability of the scheme, and also the control of the divergence error.
\Cref{sec:conclusion} concludes the paper with final remarks.

\section{Active flux scheme for the MHD}\label{sec:2d_af_scheme}
\begin{figure}[hptb!]
	\centering
	\includegraphics[width=0.35\linewidth]{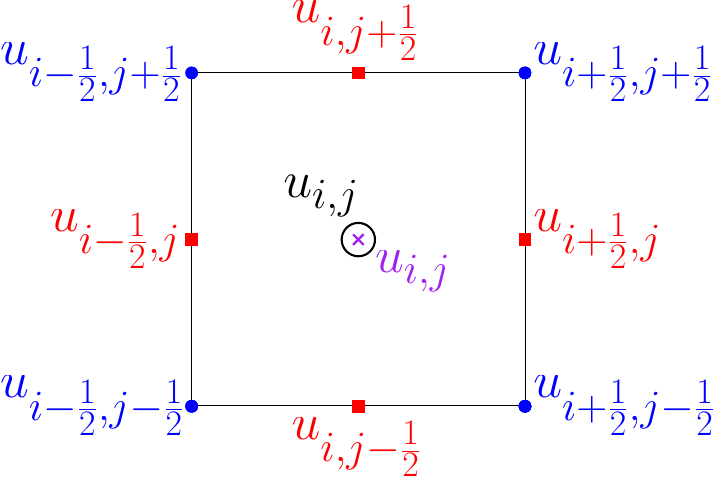}
	\caption{The DoFs for the third-order AF method: cell average (circle), face-centered values (squares), values at corners (dots). Note that the cell-centered point value $u_{i,j}$ (cross) is used in constructing the scheme, but does not belong to the DoFs.}
	\label{fig:2d_af_dofs}
\end{figure}

This section presents the 2D semi-discrete AF methods for the modified MHD system \eqref{eq:2d_mhd_gp}.
The SSP-RK3 method is used to obtain the fully-discrete scheme.
Without loss of generality, assume that a 2D computational domain is divided into $N_1\times N_2$ uniform cells,
$I_{i,j} = [x_{\xl}, x_{\xr}]\times [y_{\yl}, y_{\yr}]$ with cell centers $(x_i,y_j)=(\frac12(x_{\xl} + x_{\xr}), \frac12(y_{\yl} + y_{\yr}))$ and mesh sizes $\Delta x, \Delta y$.
The DoFs contain the cell averages and point values of the numerical solution $\bU_h(x,y,t)$, defined as
\begin{equation*}
	\overline{\bU}_{i,j}(t) = \frac{1}{\Delta x\Delta y}\int_{I_{i,j}} \bU_h(x,y,t) ~\dd \bx,\quad
	\bU_{\zeta}(t) = \bU_h(\bx_{\zeta}, t),
\end{equation*}
with $\bx=(x,y)$, $\zeta = (i+\frac12,j), (i,j+\frac12), (\xr,\yr)$.
Figure \ref{fig:2d_af_dofs} shows the locations of the DoFs for some variable $u$.
Now, let us introduce some finite difference operators, which will be used in the construction of our AF scheme.

\subsection{Finite difference operators}
Recall the quadrature points and weights of Simpson's rule in the interval $[-\frac12,\frac12]$,
\begin{equation*}
	\xi_1 = -\frac12, ~\xi_2 = 0, ~\xi_3 = \frac12,~\text{and}~
	\omega_1 = \frac16, ~\omega_2 = \frac23, ~\omega_3 = \frac16,
\end{equation*}
then the $3\times3$ tensor product quadrature points in the cell $I_{i,j}$ are
\begin{equation}\label{eq:cell_quadrature_point}
	(x_i + \xi_l\Delta x, ~ y_j + \xi_m\Delta y), \quad l,m=1,2,3.
\end{equation}
Denote the numerical solution at the $(l,m)$th quadrature points as $\bU_{i,j}^{l,m}$.
One can retrieve them easily from the DoFs since their locations coincide except for $l=m=1$, i.e., the solution at the cell center.
To obtain that, reconstruct a bi-parabolic polynomial in the cell $I_{i,j}$ using the cell average and $8$ point values on the edges \cite{Abgrall_2025_semi_JSC}, then the cell-centered solution is
\begin{align*}%\label{eq:2d_cell_center_parabolic}
	\bU_{i,j}^{1,1} = \bU_{i,j} = \frac{1}{16}\Big[
	36\overline{\bU}_{i,j} &- 4\left(\bU_{\xl,j}+\bU_{\xr,j}+\bU_{i,\yl}+\bU_{i,\yr}\right)\nonumber \\
	&- \left(\bU_{\xl,\yl} + \bU_{\xr,\yl} + \bU_{\xl,\yr} + \bU_{\xr,\yr}\right)
	\Big]. 
\end{align*}
The following finite difference operators for the first-order derivatives can be obtained by differentiating the bi-parabolic reconstruction in the cell $I_{i,j}$,
\begin{align}
	\left(D_1^{+} u\right)_{\xr,j_1} &= \dfrac{1}{\Delta x}\left( u_{\xl,j_1} - 4 u_{i,j_1} + 3 u_{\xr,j_1} \right), \nonumber\\
	\left(D_1^{-} u\right)_{\xl,j_1} &= \dfrac{1}{\Delta x}\left(- 3 u_{\xl,j_1} + 4 u_{i,j_1} - u_{\xr,j_1} \right), \nonumber\\
	\left(D_1 u\right)_{i,j_1} &= \dfrac{1}{\Delta x}\left( u_{\xr,j_1} - u_{\xl,j_1} \right), \nonumber\\
	\left(D_2^{+} u\right)_{i_1,\yr} &= \dfrac{1}{\Delta y}\left( u_{i_1,\yl} - 4 u_{i_1,j} + 3 u_{i_1,\yr} \right), \nonumber\\
	\left(D_2^{-} u\right)_{i_1,\yl} &= \dfrac{1}{\Delta y}\left(- 3 u_{i_1,\yl} + 4 u_{i_1,j} - u_{i_1,\yr} \right), \nonumber\\
	\left(D_2 u\right)_{i_1,j} &= \dfrac{1}{\Delta y}\left( u_{i_1,\yr} - u_{i_1,\yl} \right), \label{eq:2d_fd}
\end{align}
where $i_1=i-\frac12,i,i+\frac12$, $j_1=j-\frac12,j,j+\frac12$.
Note that $D^{\pm}_{\ell}, \ell=1,2$ are one-sided finite difference operators, while $D_{\ell}, \ell=1,2$ are central finite difference operators, and they are exact for bi-parabolic polynomials.

\subsection{Evolution of cell average}
The cell average is evolved following the finite volume method
\begin{equation*}%\label{eq:2d_semi_av}
	\frac{\dd \overline{\bU}_{i,j}}{\dd t} =
	-\frac{1}{\Delta x}\left(\widehat{\bF}_{1,\xr,j} - \widehat{\bF}_{1,\xl,j}\right)
	-\frac{1}{\Delta y}\left(\widehat{\bF}_{2,i,\yr} - \widehat{\bF}_{2,i,\yl}\right)
	-\bm{S}_{i,j},
\end{equation*}
where $\widehat{\bF}_{1,\xr,j}$ and $\widehat{\bF}_{2,i,\yr}$ are the numerical fluxes
\begin{equation*}%\label{eq:2d_num_fluxes}
	\widehat{\bF}_{1,\xr,j} = \frac{1}{\Delta y}\int_{y_{\yl}}^{y_{\yr}} \bF_1(\bU_h(x_{\xr},y))~\dd y,~
	\widehat{\bF}_{2,i,\yr} = \frac{1}{\Delta x}\int_{x_{\xl}}^{x_{\xr}} \bF_2(\bU_h(x,y_{\yr}))~\dd x,
\end{equation*}
and $\bm{S}$ is the discretization of the Godunov-Powell source term.
Using Simpson's rule, this paper uses the third-order numerical flux in the $x$-direction
\begin{equation*}%\label{eq:flux_simpson}
	\widehat{\bF}_{1,\xr,j}
	= \frac16\left[\bF_1(\bU_{\xr,\yl}) + 4\bF_1(\bU_{\xr,j}) + \bF_1(\bU_{\xr,\yr})\right].
\end{equation*}
We propose the following third-order discretization for the nonconservative source term
\begin{equation}\label{eq:semi_av_src}
	\bm{S}_{i,j} = \sum_{l, m=1}^{3} \omega_l\omega_m(\divB)_{i,j}^{l, m}\bm{\Psi}(\bU_{i,j}^{l, m}),
\end{equation}
where $(\divB)_{i,j}^{l,m}$ is the discrete divergence at the $(l,m)$th quadrature point defined in \eqref{eq:cell_quadrature_point}.
The derivatives in the discrete divergence are computed by the finite difference \eqref{eq:2d_fd}, i.e.,
\begin{align*}
	&(\divB)_{i,j}^{0,0}
	= \left(\frac{\partial B_1}{\partial x}\right)_{\xl,\yl} + \left(\frac{\partial B_2}{\partial y}\right)_{\xl,\yl}
	= \left(D_1^{-}B_1\right)_{\xl,\yl} +  \left(D_2^{-}B_2\right)_{\xl,\yl}, \\
	&(\divB)_{i,j}^{0,1}
	= \left(\frac{\partial B_1}{\partial x}\right)_{\xl,j} + \left(\frac{\partial B_2}{\partial y}\right)_{\xl,j}
	= \left(D_1^{-}B_1\right)_{\xl,j} +  \left(D_2B_2\right)_{\xl,j}, \\
	&(\divB)_{i,j}^{0,2}
	= \left(\frac{\partial B_1}{\partial x}\right)_{\xl,\yr} + \left(\frac{\partial B_2}{\partial y}\right)_{\xl,\yr}
	= \left(D_1^{-}B_1\right)_{\xl,\yr} +  \left(D_2^{+}B_2\right)_{\xl,\yr}, \\
    &(\divB)_{i,j}^{1,0}
	= \left(\frac{\partial B_1}{\partial x}\right)_{i,\yl} + \left(\frac{\partial B_2}{\partial y}\right)_{i,\yl}
	= \left(D_1 B_1\right)_{i,\yl} +  \left(D_2^{-}B_2\right)_{i,\yl}, \\
	&(\divB)_{i,j}^{1,1}
	= \left(\frac{\partial B_1}{\partial x}\right)_{i,j} + \left(\frac{\partial B_2}{\partial y}\right)_{i,j}
	= \left(D_1 B_1\right)_{i,j} +  \left(D_2B_2\right)_{i,j}, \\
	&(\divB)_{i,j}^{1,2}
	= \left(\frac{\partial B_1}{\partial x}\right)_{i,\yr} + \left(\frac{\partial B_2}{\partial y}\right)_{i,\yr}
	= \left(D_1 B_1\right)_{i,\yr} +  \left(D_2^{+}B_2\right)_{i,\yr}, \\
    &(\divB)_{i,j}^{2,0}
	= \left(\frac{\partial B_1}{\partial x}\right)_{\xr,\yl} + \left(\frac{\partial B_2}{\partial y}\right)_{\xr,\yl}
	= \left(D_1^{+}B_1\right)_{\xr,\yl} +  \left(D_2^{-}B_2\right)_{\xr,\yl}, \\
	&(\divB)_{i,j}^{2,1}
	= \left(\frac{\partial B_1}{\partial x}\right)_{\xr,j} + \left(\frac{\partial B_2}{\partial y}\right)_{\xr,j}
	= \left(D_1^{+}B_1\right)_{\xr,j} +  \left(D_2B_2\right)_{\xr,j}, \\
	&(\divB)_{i,j}^{2,2}
	= \left(\frac{\partial B_1}{\partial x}\right)_{\xr,\yr} + \left(\frac{\partial B_2}{\partial y}\right)_{\xr,\yr}
	= \left(D_1^{+}B_1\right)_{\xr,\yr} +  \left(D_2^{+}B_2\right)_{\xr,\yr}.
\end{align*}

\subsection{Evolution of point value}\label{sec:pnt_evolution}
This paper adopts the LLF FVS proposed in \cite{Duan_2025_Active_SJSC} for the discretization of the flux derivative, which was shown to be superior to the Jacobian splitting \cite{Abgrall_2025_semi_JSC} for strong discontinuities.
For the point value at the corner $(x_{\xr}, y_{\yr})$, the discretization is
\begin{equation*}%\label{eq:2d_semi_fvs_node}
	\dfrac{\dd \bU_{\xr,\yr}}{\dd t} 
	= - \sum_{\ell=1}^2 \left[D_\ell^+\bF_\ell^{+}(\bU)+ D_\ell^-\bF_\ell^{-}(\bU) \right]_{\xr,\yr}
	- (\divB)_{\xr,\yr}\bm{\Psi}(\bU_{\xr,\yr}),
\end{equation*}
where 
\begin{equation*}
	(\divB)_{\xr,\yr} = \sum_{\ell=1}^2 \frac12\left[ (D_\ell^{+}{B_\ell})_{\xr,\yr} + (D_\ell^{-}{B_\ell})_{\xr,\yr} \right],
\end{equation*}
and the finite difference operator is performed on each component.
Here $\bF_\ell^\pm = \frac12(\bF_\ell(\bU) \pm \alpha_\ell \bU)$ is obtained from the LLF FVS,
and the coefficient is chosen as
\begin{align*}%\label{eq:fvs_alpha}
	&(\alpha_1)_{\xr,\yr} = \max \left\{\varrho_1(\bU_{\xl,\yr}), \varrho_1(\bU_{i,\yr}), \varrho_1(\bU_{\xr,\yr}), \varrho_1(\bU_{i+1,\yr}), \varrho_1(\bU_{i+\frac32,\yr})\right\}, \\
	&(\alpha_2)_{\xr,\yr} = \max \left\{\varrho_2(\bU_{\xr,\yl}), \varrho_2(\bU_{\xr,j}), \varrho_2(\bU_{\xr,\yr}), \varrho_2(\bU_{\xr,j+1}), \varrho_2(\bU_{\xr,j+\frac32})\right\},
\end{align*}
i.e., the maximal spectral radius $\varrho_\ell$ of $\partial\bF_\ell/\partial\bU$ across the spatial stencil.
Note that we use upwind finite difference here.
One can verify that, see e.g. \cite{Powell_1994_approximate}, $\varrho_\ell = \abs{v_\ell} + c_{f,\ell}, \ell=1,2$, with
\begin{equation}\label{eq:cf}
	c_{f,\ell} = \sqrt{\frac12\left(c^2 + \frac{\norm{\bB}^2}{\rho} + \sqrt{\left(c^2 + \frac{\norm{\bB}^2}{\rho}\right)^2 - 4\frac{c^2 B_\ell^2}{\rho} } \right)},\quad
	c = \sqrt{\frac{\gamma p}{\rho}}.
\end{equation}

For the face-centered point values at $(x_{\xr}, y_{j})$ and $(x_{i}, y_{\yr})$, their discretizations are
\begin{align*}
	&\dfrac{\dd \bU_{\xr,j}}{\dd t} 
	= - \left(D_1^+\bF_1^{+} + D_1^-\bF_1^{-} \right)_{\xr,j}
	- \left(D_2\bF_2\right)_{\xr,j}
	- (\divB)_{\xr,j}\bm{\Psi}(\bU_{\xr,j}), \\
	&\dfrac{\dd \bU_{i,\yr}}{\dd t} 
	= - \left(D_1\bF_1\right)_{i,\yr}
	- \left(D_2^+\bF_2^{+}+ D_2^-\bF_2^{-} \right)_{i,\yr}
	- (\divB)_{i,\yr}\bm{\Psi}(\bU_{i,\yr}),
\end{align*}
where
\begin{align*}
	&(\divB)_{\xr,j} = \frac12\left[ (D_1^{+}{B_1})_{\xr,j} + (D_1^{-}{B_1})_{\xr,j} \right] + (D_2{B_2})_{\xr,j}, \\
	&(\divB)_{i,\yr} = (D_1{B_1})_{i,\yr} + \frac12\left[ (D_2^{+}{B_2})_{i,\yr} + (D_2^{-}{B_2})_{i,\yr} \right],
\end{align*}
and the coefficients in the FVS can be obtained similarly.

\begin{remark}
	Note that $\frac12\left( D_\ell^{+}{B_\ell} + D_\ell^{-}{B_\ell} \right)$ is a central finite difference for the derivative of $B_\ell$.
	We will show in Example \ref{ex:2d_accuracy} that such an approximation is stable, while the upwind discretization based on the velocity direction leads to instability.
\end{remark}

\begin{remark}
	The discretization of the source term for the cell average only depends on the DoFs in the cell $I_{i,j}$, and the discretization for the source term in the point value does not enlarge the spatial stencil,
	thus our scheme keeps the compactness of the original AF scheme \cite{Abgrall_2025_semi_JSC} on the Cartesian mesh.
\end{remark}

\section{Limitings for the active flux scheme}\label{sec:2d_limitings}
For the ideal MHD equations, the admissible state set is
\begin{equation*}%\label{eq:2d_mhd_g}
	\mathcal{G} = \left\{\bU = \left(\rho, \rho\bv, \bB, E\right) ~\Big|~ \rho > 0,~ p = (\gamma-1)\left(E - \frac{\norm{\rho \bv}^2}{2\rho} - \frac{\norm{\bB}^2}{2}\right) > 0 \right\},
\end{equation*}
which is convex, see e.g. \cite{Cheng_2013_Positivity_JCP}.
\begin{definition}
	An AF scheme is called \emph{positivity-preserving} (PP) if starting from admissible cell averages and point values in $\mathcal G$, the cell averages and point values stay in $\mathcal G$ at the next time step.
\end{definition}
As the DoFs in the AF scheme include cell averages and point values,
one must design suitable limitings for both of them to achieve the PP property.
The idea is to blend the high-order AF scheme with a PP scheme \cite{Duan_2025_Active_SJSC}.
To begin with, let us review the PP property of the first-order LLF scheme with the Godunov-Powell source term.

\subsection{First-order positivity-preserving LLF scheme}
Consider the following scheme for \eqref{eq:2d_mhd_gp},
\begin{align}
	\overline{\bU}_{i,j}^{\texttt{LLF}} = \overline{\bU}_{i,j}^{n}
	&- \frac{\Delta t^n}{\Delta x}\left[\widehat{\bF}_{1,\xr,j}^{\texttt{LLF}}(\overline{\bU}_{i,j}^n, \overline{\bU}_{i+1,j}^n)
	- \widehat{\bF}_{1,\xl,j}^{\texttt{LLF}}(\overline{\bU}_{i-1,j}^n, \overline{\bU}_{i,j}^n)\right] \nonumber\\
	&- \frac{\Delta t^n}{\Delta y}\left[\widehat{\bF}_{2,i,\yr}^{\texttt{LLF}}(\overline{\bU}_{i,j}^n, \overline{\bU}_{i,j+1}^n)
	- \widehat{\bF}_{2,i,\yl}^{\texttt{LLF}}(\overline{\bU}_{i,j-1}^n, \overline{\bU}_{i,j}^n)\right]
	- \Delta t^n \bm{S}_{i,j}^{\texttt{LLF}},
	\label{eq:2d_llf}
\end{align}
with the LLF flux
\begin{equation}\label{eq:2d_llf_flux}
	\widehat{\bF}_\ell^{\texttt{LLF}}(\bU, \widetilde{\bU}) = \frac12\left[\bF_\ell(\bU) + \bF_\ell(\widetilde{\bU})
	- \alpha_{\ell}^{\texttt{LLF}}(\widetilde{\bU} - \bU)\right],
\end{equation}
and the discretized source term
\begin{equation*}
	\bm{S}_{i,j}^{\texttt{LLF}} = (\nabla\cdot \overline{\bm{B}}^n)_{i,j}\bm{\Psi}(\overline{\bU}_{i,j}^n),
\end{equation*}
where the discrete divergence is computed by the central finite difference
\begin{equation*}
	(\nabla\cdot \overline{\bm{B}}^n)_{i,j} =  \left[\frac{(\overline{B_1})_{i+1,j}^n - (\overline{B_1})_{i-1,j}^n}{2\Delta x}
	+ \frac{(\overline{B_2})_{i,j+1}^n - (\overline{B_2})_{i,j-1}^n}{2\Delta y}\right].
\end{equation*}

\begin{lemma}[\cite{Wu_2019_Provably_NM}]\label{lem:llf}
	Assume that the parameters in \eqref{eq:2d_llf_flux} satisfy for $\ell=1,2$,
	\begin{equation*}
		\alpha_\ell^{\texttt{LLF}} \geqslant \max\{\varrho_\ell(\bU), \varrho_\ell(\widetilde{\bU}), \alpha_{\ell,\ast}(\bU, \widetilde{\bU}),  \alpha_{\ell,\ast}(\widetilde{\bU}, \bU)\},
	\end{equation*}
	with
	\begin{equation*}
		\alpha_{\ell,\ast}(\bU, \widetilde{\bU}) = \max\left\{\abs{v_\ell},
		\frac{\abs{\sqrt{\rho}v_\ell + \sqrt{\widetilde{\rho}}\widetilde{v_\ell}}}{\sqrt{\rho} + \sqrt{\widetilde{\rho}}} \right\}
		+ \max\{c_{f,\ell}, \widetilde{c_{f,\ell}}\}
		+ \frac{\norm{\bB - \widetilde{\bB}}}{\sqrt{\rho} + \sqrt{\widetilde{\rho}}},
	\end{equation*}
	where $c_{f,\ell}$ and $\widetilde{c_{f,\ell}}$ are evaluated based on $\bU$ and $\widetilde{\bU}$ according to \eqref{eq:cf}, respectively.
	Given $\overline{\bU}_{i,j}^n\in\mathcal{G}, \forall i,j$, then the solution of the first-order scheme \eqref{eq:2d_llf} is PP, i.e., $\overline{\bU}_{i,j}^{\texttt{LLF}}\in\mathcal{G}$ under the CFL condition
	\begin{equation*}
		\Delta t^n\left( \frac{\alpha_{1,\xl,j}^{\texttt{LLF}}}{\Delta x}
		+ \frac{\alpha_{1,\xr,j}^{\texttt{LLF}}}{\Delta x}
		+ \frac{\alpha_{2,i,\yl}^{\texttt{LLF}}}{\Delta y}
		+ \frac{\alpha_{2,i,\yr}^{\texttt{LLF}}}{\Delta y}
		+ \frac{ \abs{(\nabla\cdot \overline{\bm{B}}^n)_{i,j}} }{\sqrt{ \bar{\rho}_{i,j}^n}} \right) \leqslant 1, ~\forall i,j.
	\end{equation*}	
\end{lemma}

In our limitings, the high-order AF scheme will be blended with the above PP LLF scheme.
If the high-order scheme equipped with the forward Euler scheme is PP,
then the high-order scheme using the SSP-RK3 is also PP since the SSP-RK3 is a convex combination of forward-Euler stages.
Thus, only the forward Euler scheme is considered below.
Note that to avoid the effect of the round-off error, we need to choose desired lower bounds for the density and pressure,
denoted by $\varepsilon_{\rho}, \varepsilon_{p}$ to be defined later,
such that $\rho\geqslant\varepsilon_{\rho}$, $p\geqslant\varepsilon_{p}$.

\subsection{Parametrized flux limiter for cell average}\label{sec:2d_limiting_average}
This section presents a flux limiting approach to enforce the PP property of the cell average update by constraining individual numerical fluxes \cite{Xu_2014_Parametrized_MC,Christlieb_2015_Positivity_SJSC,Christlieb_2016_High_JCP}.
Let $\varepsilon_\rho = \min\{10^{-13}, \rho(\overline{\bU}_{i,j}^{\texttt{LLF}})\}$, $\varepsilon_p = \min\{10^{-13}, p(\overline{\bU}_{i,j}^{\texttt{LLF}})\}$.
As the first-order solution satisfies $\rho(\overline{\bU}_{i,j}^{\texttt{LLF}}) \geqslant \varepsilon_\rho$, $p(\overline{\bU}_{i,j}^{\texttt{LLF}}) \geqslant \varepsilon_p$, one can find limited fluxes and source term by blending the high-order and first-order parts as
\begin{align*}
	\widehat{\bF}_{1,i\pm\frac12,j}^{\texttt{Lim}} &= \theta_{i\pm\frac12,j} \widehat{\bF}_{1,i\pm\frac12,j} + (1-\theta_{i\pm\frac12,j})\widehat{\bF}_{1,i\pm\frac12,j}^{\texttt{LLF}}, \\
	\widehat{\bF}_{2,i,j\pm\frac12}^{\texttt{Lim}} &= \theta_{i,j\pm\frac12} \widehat{\bF}_{2,i,j\pm\frac12} + (1-\theta_{i,j\pm\frac12})\widehat{\bF}_{2,i,j\pm\frac12}^{\texttt{LLF}}, \\
	\bm{S}_{i,j}^{\texttt{Lim}} &= \theta_{i,j} \bm{S}_{i,j} + (1-\theta_{i,j})\bm{S}_{i,j}^{\texttt{LLF}},
\end{align*}
such that the solution of the following limited scheme
\begin{equation}\label{eq:2d_av_limited}
	\overline{\bU}_{i,j}^{\texttt{Lim}} = \overline{\bU}_{i,j}^{n}
	-\frac{\Delta t^n}{\Delta x}\left(\widehat{\bF}_{1,\xr,j}^{\texttt{Lim}} - \widehat{\bF}_{1,\xl,j}^{\texttt{Lim}}\right)
	-\frac{\Delta t^n}{\Delta y}\left(\widehat{\bF}_{2,i,\yr}^{\texttt{Lim}} - \widehat{\bF}_{2,i,\yl}^{\texttt{Lim}}\right)
	-\Delta t^n\bm{S}_{i,j}^{\texttt{Lim}}
\end{equation}
satisfies $\rho(\overline{\bU}_{i,j}^{\texttt{Lim}})\geqslant \varepsilon_\rho$, $p(\overline{\bU}_{i,j}^{\texttt{Lim}})\geqslant \varepsilon_p$.
The coefficients $\theta_{i\pm\frac12,j}, \theta_{i,j\pm\frac12}, \theta_{i,j}$ should stay in $[0,1]$, and be as close to $1$ as possible, so that the high-order terms are used as much as possible to maintain the accuracy.
Our PP limiting consists of the following three steps,
and mainly follows the notations in \cite{Christlieb_2015_Positivity_SJSC,Christlieb_2016_High_JCP}.
Note that their procedure is two steps, as the source term is not included.

\noindent(1) \textbf{Limit the source term} such that
\begin{align*}
	\overline{\bU}_{i,j}^{\texttt{Lim,1}}
	&= \overline{\bU}_{i,j}^{n}
	-\frac{\Delta t^n}{\Delta x}\left(\widehat{\bF}_{1,\xr,j}^{\texttt{LLF}} - \widehat{\bF}_{1,\xl,j}^{\texttt{LLF}}\right)
	-\frac{\Delta t^n}{\Delta y}\left(\widehat{\bF}_{2,i,\yr}^{\texttt{LLF}} - \widehat{\bF}_{2,i,\yl}^{\texttt{LLF}}\right)
	-\Delta t^n\bm{S}_{i,j}^{\texttt{Lim}} \\
	&= \theta_{i,j}\overline{\bU}_{i,j}^{\texttt{src}}
	+ (1-\theta_{i,j})\overline{\bU}_{i,j}^{\texttt{LLF}}
\end{align*}
is PP, where
\begin{equation*}
	\overline{\bU}_{i,j}^{\texttt{src}} = \overline{\bU}_{i,j}^{n}
	-\frac{\Delta t^n}{\Delta x}\left(\widehat{\bF}_{1,\xr,j}^{\texttt{LLF}} - \widehat{\bF}_{1,\xl,j}^{\texttt{LLF}}\right)
	-\frac{\Delta t^n}{\Delta y}\left(\widehat{\bF}_{2,i,\yr}^{\texttt{LLF}} - \widehat{\bF}_{2,i,\yl}^{\texttt{LLF}}\right)
	-\Delta t^n\bm{S}_{i,j}.
\end{equation*}
Since the first component of the source term is zero as seen in \eqref{eq:gp_src},
the density of $\overline{\bU}_{i,j}^{\texttt{Lim,1}}$ or $\overline{\bU}_{i,j}^{\texttt{src}}$ is the same as $\overline{\bU}_{i,j}^{\texttt{LLF}}$, which automatically satisfies $\rho(\overline{\bU}_{i,j}^{\texttt{Lim,1}}) \geqslant \varepsilon_\rho$.
Due to the concavity of the pressure \cite{Cheng_2013_Positivity_JCP}, one has
\begin{equation*}
	p(\overline{\bU}_{i,j}^{\texttt{Lim,1}}) \geqslant \theta_{i,j}p(\overline{\bU}_{i,j}^{\texttt{src}}) + (1-\theta_{i,j})p(\overline{\bU}_{i,j}^{\texttt{LLF}}).
\end{equation*}
Define
\begin{equation*}
	\theta_{i,j} = \min\left\{1, \frac{p(\overline{\bU}_{i,j}^{\texttt{LLF}}) - \varepsilon_p}
	{p(\overline{\bU}_{i,j}^{\texttt{LLF}}) - p(\overline{\bU}_{i,j}^{\texttt{src}})} \right\},
\end{equation*}
then it is easy to verify that $p(\overline{\bU}_{i,j}^{\texttt{Lim,1}}) \geqslant \varepsilon_p$, thus $\overline{\bU}_{i,j}^{\texttt{Lim,1}} \in \mathcal{G}$.

\noindent(2) \textbf{Find candidate parameters} $\Lambda_{I_{i,j}, \texttt{L}}$, $\Lambda_{I_{i,j}, \texttt{R}}$, $\Lambda_{I_{i,j}, \texttt{D}}$, $\Lambda_{I_{i,j}, \texttt{U}}$ as close to $1$ as possible in each cell $I_{i,j}$ such that for all
\begin{equation*}
	(\theta_{\texttt{L}}, \theta_{\texttt{R}}, \theta_{\texttt{D}}, \theta_{\texttt{U}}) \in [0, \Lambda_{I_{i,j}, \texttt{L}}]\times[0, \Lambda_{I_{i,j}, \texttt{R}}]\times [0, \Lambda_{I_{i,j}, \texttt{D}}]\times[0, \Lambda_{I_{i,j}, \texttt{U}}],
\end{equation*}
the limited solution
\begin{equation*}
	\overline{\bU}_{i,j}^{\texttt{Lim,2}}(\theta_{\texttt{L}}, \theta_{ \texttt{R}}, \theta_{\texttt{D}}, \theta_{\texttt{U}})
	= \overline{\bU}_{i,j}^{\texttt{Lim},1}
	+ \theta_{\texttt{L}}\bm{\mathcal{H}}_\texttt{L}
	+ \theta_{\texttt{R}}\bm{\mathcal{H}}_\texttt{R}
	+ \theta_{\texttt{D}}\bm{\mathcal{H}}_\texttt{D}
	+ \theta_{\texttt{U}}\bm{\mathcal{H}}_\texttt{U}
\end{equation*}
is PP, where the anti-diffusive fluxes are given by
\begin{align*}
	&\bm{\mathcal{H}}_\texttt{L} = 
	\frac{\Delta t^n}{\Delta x}\left(\widehat{\bF}_{1,\xl,j} - \widehat{\bF}_{1,\xl,j}^{\texttt{LLF}}\right),~
	\bm{\mathcal{H}}_\texttt{R} = 
	- \frac{\Delta t^n}{\Delta x}\left(\widehat{\bF}_{1,\xr,j} - \widehat{\bF}_{1,\xr,j}^{\texttt{LLF}}\right), \nonumber\\
	&\bm{\mathcal{H}}_\texttt{D} = 
	\frac{\Delta t^n}{\Delta y}\left(\widehat{\bF}_{2,i,\yl} - \widehat{\bF}_{2,i,\yl}^{\texttt{LLF}}\right),~
	\bm{\mathcal{H}}_\texttt{U} = 
	- \frac{\Delta t^n}{\Delta y}\left(\widehat{\bF}_{2,i,\yr} - \widehat{\bF}_{2,i,\yr}^{\texttt{LLF}}\right).
\end{align*}
Because the following two sets
\begin{equation*}
	S_\rho = \{ (\theta_{\texttt{L}}, \theta_{ \texttt{R}}, \theta_{\texttt{D}}, \theta_{\texttt{U}}) \in[0,1]^4 ~|~ \rho(\overline{\bU}_{i,j}^{\texttt{Lim,2}}(\theta_{\texttt{L}}, \theta_{ \texttt{R}}, \theta_{\texttt{D}}, \theta_{\texttt{U}})) \geqslant \varepsilon_\rho \}
\end{equation*}
and
\begin{equation*}
	S_p = \{ (\theta_{\texttt{L}}, \theta_{ \texttt{R}}, \theta_{\texttt{D}}, \theta_{\texttt{U}}) \in[0,1]^4 ~|~ \rho(\overline{\bU}_{i,j}^{\texttt{Lim,2}}(\theta_{\texttt{L}}, \theta_{ \texttt{R}}, \theta_{\texttt{D}}, \theta_{\texttt{U}})) \geqslant \varepsilon_\rho ~\text{and}~
	p(\overline{\bU}_{i,j}^{\texttt{Lim,2}}(\theta_{\texttt{L}}, \theta_{ \texttt{R}}, \theta_{\texttt{D}}, \theta_{\texttt{U}})) \geqslant \varepsilon_p \}
\end{equation*}
are convex \cite{Christlieb_2015_Positivity_SJSC,Christlieb_2016_High_JCP},
one can determine the parameters $\Lambda_{I_{i,j}, \texttt{L}}$, $\Lambda_{I_{i,j}, \texttt{R}}$, $\Lambda_{I_{i,j}, \texttt{D}}$, $\Lambda_{I_{i,j}, \texttt{U}}$ in the following two steps.

\noindent $\bullet$ Find a rectangular subset $R_\rho = [0, \Lambda_{\texttt{L}}^\rho]\times[0, \Lambda_{\texttt{R}}^\rho]\times[0, \Lambda_{\texttt{D}}^\rho]\times[0, \Lambda_{\texttt{U}}^\rho]$ of $S_\rho$.
To be specific,
\begin{equation*}
	\Lambda_{\texttt{I}}^\rho = \begin{cases}
	\min\left\{ 1,~ \dfrac{\rho(\overline{\bU}_{i,j}^{\texttt{Lim},1}) - \varepsilon_\rho}{10^{-12} - \sum\limits_{\texttt{J}, \rho(\bm{\mathcal{H}}_\texttt{J}) < 0}\rho(\bm{\mathcal{H}}_\texttt{J})} \right\}, &\text{if}~\rho(\bm{\mathcal{H}}_\texttt{I}) < 0, \\
	1, &\text{otherwise},
	\end{cases}
\end{equation*}
where $\texttt{I}$ and $\texttt{J}$ take values in $\texttt{L},\texttt{R},\texttt{D},\texttt{U}$.

\noindent $\bullet$ Shrink the rectangle $R_\rho$ to make it stay within $S_p$.
Let the vertices of $R_\rho$ be $\bm{A}^{k_{\texttt{L}}, k_{\texttt{R}}, k_{\texttt{D}}, k_{\texttt{U}}}$ with $k_{\texttt{I}} = 0$ or $1$, such that the $\texttt{I}$th component of $\bm{A}^{k_{\texttt{L}}, k_{\texttt{R}}, k_{\texttt{D}}, k_{\texttt{U}}}$ is $\Lambda_{\texttt{I}}^\rho$ for $k_{\texttt{I}} = 1$ otherwise $0$.
For each $(k_{\texttt{L}}, k_{\texttt{R}}, k_{\texttt{D}}, k_{\texttt{U}})$, if $p(\bm{A}^{k_{\texttt{L}}, k_{\texttt{R}}, k_{\texttt{D}}, k_{\texttt{U}}}) \geqslant \varepsilon_p$, set the new vertex as $\bm{B}^{k_{\texttt{L}}, k_{\texttt{R}}, k_{\texttt{D}}, k_{\texttt{U}}} = \bm{A}^{k_{\texttt{L}}, k_{\texttt{R}}, k_{\texttt{D}}, k_{\texttt{U}}}$.
Otherwise, by solving a cubic equation to get the smallest positive value $r$ satisfying $p(r\bm{A}^{k_{\texttt{L}}, k_{\texttt{R}}, k_{\texttt{D}}, k_{\texttt{U}}}) \geqslant \varepsilon_p$, set the new vertex as $\bm{B}^{k_{\texttt{L}}, k_{\texttt{R}}, k_{\texttt{D}}, k_{\texttt{U}}} = r\bm{A}^{k_{\texttt{L}}, k_{\texttt{R}}, k_{\texttt{D}}, k_{\texttt{U}}}$.
Here we use the Newton method which converges within $4$ iterations in the numerical tests.
Finally, let us find a rectangular subset inside the convex polygon with vertices $\bm{B}^{k_{\texttt{L}}, k_{\texttt{R}}, k_{\texttt{D}}, k_{\texttt{U}}}$ by
\begin{equation*}
	\Lambda_{I_{i,j}, \texttt{I}} = \min\limits_{(k_{\texttt{L}}, k_{\texttt{R}}, k_{\texttt{D}}, k_{\texttt{U}}), k_{\texttt{I}} = 1}
	B_{\texttt{I}}^{k_{\texttt{L}}, k_{\texttt{R}}, k_{\texttt{D}}, k_{\texttt{U}}},
\end{equation*}
where $B_{\texttt{I}}^{k_{\texttt{L}}, k_{\texttt{R}}, k_{\texttt{D}}, k_{\texttt{U}}}$ denotes the $\texttt{I}$th component of $\bm{B}^{k_{\texttt{L}}, k_{\texttt{R}}, k_{\texttt{D}}, k_{\texttt{U}}}$.

\noindent(3) \textbf{Compute the unique blending coefficients} at cell interfaces by
\begin{equation*}
	\theta_{\xr,j} = \min\{ \Lambda_{I_{i,j}, \texttt{R}}, \Lambda_{I_{i+1,j}, \texttt{L}} \},\quad
	\theta_{i,\yr} = \min\{ \Lambda_{I_{i,j}, \texttt{U}}, \Lambda_{I_{i,j+1}, \texttt{D}} \}.
\end{equation*}

\begin{remark}\rm
	The limited scheme \eqref{eq:2d_av_limited} keeps mass conservation as the fluxes at the cell interfaces are unique.
	The conservation of momentum or total energy is not guaranteed due to the Godunov-Powell source term.
\end{remark}

\subsection{Scaling limiter for point value}\label{sec:2d_limiting_point}
The PP limiting for the point value is borrowed from \cite{Duan_2025_Active_SJSC},
i.e., blending the whole state of conservative variables directly by using the simple scaling limiter \cite{Liu_1996_Nonoscillatory_SJNA},
as there is no conservation requirement on the point value update.

The first step is to define suitable first-order LLF schemes.
For the point value at the corner, one can choose
\begin{align*}%\label{eq:2d_llf_node}
	\bU_{\xr,\yr}^{\texttt{LLF}} =&\ \bU_{\xr,\yr}^{n}
	- \dfrac{\Delta t^n}{\Delta x}
	\left(\widehat{\bF}^{\texttt{LLF}}_{1,i+1,\yr}
	- \widehat{\bF}^{\texttt{LLF}}_{1,i,\yr}\right)
	- \dfrac{\Delta t^n}{\Delta y}
	\left(\widehat{\bF}^{\texttt{LLF}}_{2,\xr,j+1}
	- \widehat{\bF}^{\texttt{LLF}}_{2,\xr,j}\right) \nonumber\\
	&- \Delta t^n\left[ \frac{(B_1)_{i+\frac32,\yr} - (B_1)_{\xl,\yr}}{2\Delta x} 
	+ \frac{(B_2)_{\xr,j+\frac32} - (B_2)_{\xr,\yl}}{2\Delta y} \right]\bm{\Psi}(\bU_{\xr,\yr}^{n}),
\end{align*}
with the LLF numerical fluxes
\begin{equation*}
	\widehat{\bF}^{\texttt{LLF}}_{1,i+1,\yr}:=\widehat{\bF}_1^{\texttt{LLF}}(\bU_{\xr,\yr}^n, \bU_{i+\frac32,\yr}^n),~
	\widehat{\bF}^{\texttt{LLF}}_{2,\xr,j+1}:=\widehat{\bF}_2^{\texttt{LLF}}(\bU_{\xr,\yr}^n, \bU_{\xr,j+\frac32}^n),
\end{equation*}
which are defined in \eqref{eq:2d_llf_flux}.

For the vertical face-centered point value, we choose the first-order LLF scheme as
\begin{align*}%\label{eq:2d_llf_facex}
	\bU_{\xr,j}^{\texttt{LLF}} =&\  \bU_{\xr,j}^{n}
	- \dfrac{\Delta t^n}{\Delta x}
	\left(\widehat{\bF}^{\texttt{LLF}}_{1,i+1,j}
	- \widehat{\bF}^{\texttt{LLF}}_{1,i,j}\right) 
	- \dfrac{\Delta t^n}{\Delta y}
	\left(\widehat{\bF}^{\texttt{LLF}}_{2,\xr,j+\frac12}
	- \widehat{\bF}^{\texttt{LLF}}_{2,\xr,j-\frac12}\right) \\
	&- \Delta t^n\left[ \frac{(B_1)_{i+\frac32,j} - (B_1)_{\xl,j}}{2\Delta x} 
	+ \frac{(B_2)_{\xr,j+\frac12} - (B_2)_{\xr,\yl}}{2\Delta y} \right]\bm{\Psi}(\bU_{\xr,j}^{n}),
\end{align*}
with the LLF numerical fluxes
\begin{equation*}
	\widehat{\bF}^{\texttt{LLF}}_{1,i+1,j}:=\widehat{\bF}_1^{\texttt{LLF}}(\bU_{\xr,j}^n, \bU_{i+\frac32,j}^n),~
	\widehat{\bF}^{\texttt{LLF}}_{2,\xr,j+\frac12}:=\widehat{\bF}_2^{\texttt{LLF}}(\bU_{\xr,j}^n, \bU_{\xr,\yr}^n).
\end{equation*}
The LLF scheme for the face-centered value on the horizontal face can be chosen as
\begin{align*}%\label{eq:2d_llf_facey}
	\bU_{i,\yr}^{\texttt{LLF}} =&\ \bU_{i,\yr}^{n}
	- \dfrac{\Delta t^n}{\Delta x}
	\left(\widehat{\bF}^{\texttt{LLF}}_{1,i+\frac12,\yr}
	- \widehat{\bF}^{\texttt{LLF}}_{1,i-\frac12,\yr}\right)
	- \dfrac{\Delta t^n}{\Delta y}
	\left(\widehat{\bF}^{\texttt{LLF}}_{2,i,j+1}
	- \widehat{\bF}^{\texttt{LLF}}_{2,i,j}\right) \\
	&- \Delta t^n\left[ \frac{(B_1)_{i+\frac12,\yr} - (B_1)_{\xl,\yr}}{2\Delta x} 
	+ \frac{(B_2)_{i,j+\frac32} - (B_2)_{i,\yl}}{2\Delta y} \right]\bm{\Psi}(\bU_{i,\yr}^{n}),
\end{align*}
with the LLF numerical fluxes
\begin{equation*}
	\widehat{\bF}^{\texttt{LLF}}_{1,\xr,\yr}:=\widehat{\bF}_1^{\texttt{LLF}}(\bU_{i,\yr}^n, \bU_{\xr,\yr}^n),~
	\widehat{\bF}^{\texttt{LLF}}_{2,i,j+1}:=\widehat{\bF}_2^{\texttt{LLF}}(\bU_{i,\yr}^n, \bU_{i,j+\frac32}^n).	
\end{equation*}
The above three first-order LLF schemes for the point values are PP according to Lemma \ref{lem:llf}.
Next, we present the PP limitings for the point value by blending the high-order AF scheme using the forward Euler stage and the LLF schemes as
\begin{equation*}
	\bU_{\sigma}^{\texttt{Lim}} = \theta_{\sigma} \bU_{\sigma}^{\texttt{H}} + (1-\theta_{\sigma}) \bU_{\sigma}^{\texttt{LLF}},
\end{equation*}
such that $\rho(\bU_{\sigma}^{\texttt{Lim}}) \geqslant  \varepsilon_\rho$, $p(\bU_{\sigma}^{\texttt{Lim}}) \geqslant  \varepsilon_p$,
where $\sigma$ denotes the locations of the point value, i.e.,
$(\xr,\yr), (\xr,j), (i,\yr)$,
and $\bU_{\sigma}^{\texttt{H}}$ is the high-order AF solution.
The lower bounds are chosen as $\varepsilon_\rho = \min\{10^{-13}, \rho(\bU_{\sigma}^{\texttt{LLF}})\}$, $\varepsilon_p = \min\{10^{-13}, p(\bU_{\sigma}^{\texttt{LLF}})\}$.

\noindent(1) \textbf{Enforce density positivity}.
%First, the high-order solution $\bU_{\sigma}^{\texttt{H}}$ is modified as $\bU_{\sigma}^{\texttt{Lim}, *}$,
%such that $\rho(\bU_{\sigma}^{\texttt{Lim}, *}) \geqslant \varepsilon^\rho$.
Choose the parameter
\begin{equation*}
	\theta^{*}_{\sigma} = \begin{cases}
		\dfrac{\rho(\bU_{\sigma}^{\texttt{LLF}})-\varepsilon_\rho}{\rho(\bU_{\sigma}^{\texttt{L}})-\rho(\bU_{\sigma}^{\texttt{H}})},&\text{if}~~ \rho(\bU_{\sigma}^{\texttt{H}}) < \varepsilon_\rho, \\
		1, &\text{otherwise}, \\
	\end{cases}
\end{equation*}
then the density component of the limited solution is modified as $\rho(\bU_{\sigma}^{\texttt{Lim}, *}) = \theta^{*}_{\sigma} \rho(\bU_{\sigma}^{\texttt{H}}) + (1-\theta^{*}_{\sigma}) \rho(\bU_{\xr}^{\texttt{LLF}}) \geqslant \varepsilon_\rho $,
with the other components remaining the same as $\bU_{\sigma}^{\texttt{H}}$.

\noindent(2) \textbf{Enforce pressure positivity}.
Modify the solution $\bU_{\sigma}^{\texttt{Lim}, *}$ as $\bU_{\sigma}^{\texttt{Lim}}$,
such that $p(\bU_{\sigma}^{\texttt{Lim}}) \geqslant \varepsilon_p$.
Let the limited solution be
\begin{equation*}%\label{eq:2d_pnt_limited_state_Euler}
	\bU_{\sigma}^{\texttt{Lim}} = \theta^{**}_{\sigma} \bU_{\sigma}^{\texttt{Lim}, *} + \left(1-\theta^{**}_{\sigma}\right) \bU_{\sigma}^{\texttt{LLF}}.
\end{equation*}
Using the concavity of pressure, we can choose the parameter as
\begin{equation*}
	\theta^{**}_{\sigma} = \begin{cases}
		\dfrac{p(\bU_{\sigma}^{\texttt{LLF}}) - \varepsilon_p}{p(\bU_{\sigma}^{\texttt{LLF}}) - p(\bU_{\sigma}^{\texttt{Lim}, *})},&\text{if}~~ p(\bU_{\sigma}^{\texttt{Lim}, *}) < \varepsilon_p, \\
		1, &\text{otherwise}. \\
	\end{cases}
\end{equation*}

\begin{remark}\label{rmk:cell_center}
To compute the high-order FVS-based point value update, we should limit the cell-centered value $\bU_{i,j}$ at the beginning of each Runge-Kutta stage.
For example, we modify $\bU_{i,j}$ as $\bU_{i,j}^{\texttt{Lim}} = \tilde{\theta}_{i,j}\bU_{i,j} + (1-\tilde{\theta}_{i,j})\overline{\bU}_{i,j}$ such that
\begin{equation*}
	\rho(\bU_{i,j}^{\texttt{Lim}}) \geqslant \min\{10^{-13}, \rho(\overline{\bU}_{i,j})\},~
	p(\bU_{i,j}^{\texttt{Lim}}) \geqslant \min\{10^{-13}, p(\overline{\bU}_{i,j})\}.
\end{equation*}
The computation of $\tilde{\theta}_{i,j}$ is similar to the procedure in this section.
\end{remark}

%\begin{proposition}
%	If the initial numerical solution $\overline{\bU}_{i,j}^0, \bU_{\sigma}^0\in\mathcal{G}$ for all $i,j,\sigma$,
%	and the time step size satisfies \eqref{eq:2d_convex_combination_dt} and \eqref{eq:2d_pnt_llf_dt},
%	then the AF methods \eqref{eq:semi_av_2d}-\eqref{eq:semi_pnt_2d} equipped with the SSP-RK3 \eqref{eq:ssp_rk3} and the BP limitings
%	preserve positive density and pressure for the Euler equations.
%\end{proposition}

%\begin{remark}
%	For uniform meshes, and if taking the maximal spectral radius of $\partial\bF_1/\partial\bU$ and $\partial\bF_2/\partial\bU$ in the domain as $\norm{\varrho_1}_\infty$ and $\norm{\varrho_2}_\infty$,
%	the following CFL condition
%	\begin{equation*}%\label{eq:2d_uniform_mesh_dt}
%		\Delta t^n \leqslant \frac14\min\left\{\dfrac{\Delta x}{\norm{\varrho_1}_\infty},
%		\dfrac{\Delta y}{\norm{\varrho_2}_\infty}
%		\right\}
%	\end{equation*}
%	fulfills the time step size constraints \eqref{eq:2d_convex_combination_dt} and \eqref{eq:2d_pnt_llf_dt}.
%\end{remark}

\subsection{Shock sensor-based limiting}\label{sec:2d_limiting_shock_sensor}
The PP limitings are not enough to suppress spurious oscillations, especially near strong shock waves.
Based on the numerical results in \cite{Duan_2025_Active_SJSC}, we found that the shock sensor-based limiting can effectively reduce oscillations.
Here we present a new shock sensor for the MHD system.
We first consider
\begin{equation*}
	(\varphi_1)_{i,j} = \dfrac{\abs{(\bar{p}_t)_{i+1,j} - 2(\bar{p}_t)_{i,j} + (\bar{p}_t)_{i-1,j}}}{\abs{(\bar{p}_t)_{i+1,j} + 2(\bar{p}_t)_{i,j} + (\bar{p}_t)_{i-1,j}}}
\end{equation*}
by replacing the fluid pressure as the total pressure $p_t = p + p_m$ in the Jameson's shock sensor \cite{Jameson_1981_Solutions_AJ}.
The second is the following modified Ducros' shock sensor \cite{Ducros_1999_Large_JCP}
\begin{equation*}
	(\varphi_2)_{i,j} = \max\left\{\dfrac{-(\nabla\cdot\bar{\bv})_{i,j}}{\sqrt{(\nabla\cdot\bar{\bv})_{i,j}^2 + (\nabla\times\bar{\bv})_{i,j}^2 + 10^{-40}}}, ~0\right\}
\end{equation*}
already used in \cite{Duan_2025_Active_SJSC}, where
\begin{align*}
	&(\nabla\cdot\bar{\bv})_{i,j}\approx \dfrac{(\bar{v}_1)_{i+1,j} - (\bar{v}_1)_{i-1,j}}{2\Delta x}
	+ \dfrac{(\bar{v}_2)_{i,j+1} - (\bar{v}_2)_{i,j-1}}{2\Delta y},\\
	&(\nabla\times\bar{\bv})_{i,j}\approx \dfrac{(\bar{v}_2)_{i+1,j} - (\bar{v}_2)_{i-1,j}}{2\Delta x}
	- \dfrac{(\bar{v}_1)_{i,j+1} - (\bar{v}_1)_{i,j-1}}{2\Delta y}.
\end{align*}
Here, $(\varphi_2)_{i,j}$ is only activated when the velocity divergence is negative.
For the MHD system, we propose to include the following discrete divergence
\begin{equation*}
	(\varphi_3)_{i,j} = \dfrac{\abs{(\bar{B}_1)_{i+1,j} - (\bar{B}_1)_{i-1,j}
		+ (\bar{B}_2)_{i,j+1} - (\bar{B}_2)_{i,j-1}}}{\abs{(\bar{B}_1)_{i,j} + (\bar{B}_2)_{i,j}} + 10^{-40}},
\end{equation*}
which takes into account the divergence error.
Note that the quantities $\bar{a}_{i,j}$ used above are recovered from the cell average $\overline{\bU}_{i,j}$.
The final blending coefficient is designed as
\begin{align*}
	&\theta_{\xr,j}^{s} = \exp\left(-\kappa \left[(\varphi_1)_{\xr,j} (\varphi_2)_{\xr,j} + (\varphi_3)_{\xr,j}\right]\right)\in (0, 1],\\
	&(\varphi_s)_{\xr,j} = \max\left\{(\varphi_s)_{i,j}, (\varphi_s)_{i+1,j}\right\}, ~s=1,2,3,
\end{align*}
where the parameter $\kappa$ determines the limiting strength.
The limited numerical flux is
\begin{equation*}%\label{eq:2d_flux_limited_mhd_ss}
	\widehat{\bF}_{1,\xr,j}^{\texttt{Lim}} = (1-\theta_{\xr,j}^{s})\widehat{\bF}_{1,\xr,j}^{\texttt{LLF}} + \theta_{\xr,j}^{s}\widehat{\bF}_{1,\xr,j}.
\end{equation*}
The high-order discretization of the Godunov-Powell source term may also introduce oscillations, thus, we choose to limit the source term as
\begin{equation*}%\label{eq:2d_src_limited_mhd_ss}
	\bm{S}_{i,j}^{\texttt{Lim}} = (1-\theta_{i,j}^{s})\bm{S}_{i,j}^{\texttt{LLF}} + \theta_{i,j}^{s}\bm{S}_{i,j},
\end{equation*}
where
\begin{equation*}
	\theta_{i,j}^{s} = \min\{ \theta_{\xl,j}^{s}, \theta_{\xr,j}^{s}, \theta_{i,\yl}^{s}, \theta_{i,\yr}^{s} \}.
\end{equation*}
Note that the shock sensor-based limiting is applied before the PP limitings.

\section{Numerical results}\label{sec:results}
This section conducts some numerical tests to verify the accuracy, PP property, and shock-capturing ability of our AF scheme.
The 1D tests are computed by using the 1D AF scheme based on the LLF FVS without the Godunov-Powell source term, since the divergence-free condition holds automatically.
More details on the 1D AF scheme can also be found in \cite{Duan_2025_Active_SJSC}.
The 2D visualization is based on a refined mesh with half the mesh size, where the values at the grid points are the cell averages or point values on the original mesh.

\begin{example}[1D Riemann problems]\label{ex:1d_rp}
	The computational domain is $[0,1]$, and two cases are considered.
	The initial data of the first case \cite{Ryu_1995_Numerical_AJ} are
	\begin{equation*}
		(\rho, \bv, \bB, p) =
		\begin{cases}
			(1, ~0, ~0, ~0, ~0.75, ~1, ~0, ~1), &\text{if}~x < 0.5, \\
			(0.125, ~0, ~0, ~0, ~0.75, ~-1, ~0, ~0.1), &\text{otherwise},
		\end{cases}
	\end{equation*}
	and the adiabatic index is $\gamma=2$.
	The initial data of the second case \cite{Brio_1988_upwind_JCP} are
	\begin{equation*}
		(\rho, \bv, \bB, p) =
		\begin{cases}
			(1.08, ~1.2, ~0.01, ~0.5, ~{2}/{\sqrt{4\pi}}, ~{3.6}/{\sqrt{4\pi}}, ~{2}/{\sqrt{4\pi}}, ~0.95), &\text{if}~x < 0.5, \\
			(1, ~0, ~0, ~0, ~{2}/{\sqrt{4\pi}}, ~{4}/{\sqrt{4\pi}}, ~{2}/{\sqrt{4\pi}}, ~1), &\text{otherwise},
		\end{cases}
	\end{equation*}
	with the adiabatic index $5/3$.
	The final time is $T=0.2$ for both cases.
	The reference solution is obtained by a second-order HLLD finite volume scheme \cite{Miyoshi_2005_multi_JCP} on a fine mesh with $5000$ cells.
	
	The results obtained with our AF scheme and $800$ cells are shown in Figures \ref{fig:1d_rp1}-\ref{fig:1d_rp2}.
	The parameter $\kappa$ in the shock sensor is chosen as $10$ and $50$ for the two cases.
	One can observe that our AF scheme can capture the discontinuities with high resolution and only a few overshoots or undershoots.
	
	\begin{figure}[htb!]
		\begin{subfigure}[b]{0.32\textwidth}
			\centering
			\includegraphics[width=1.0\linewidth]{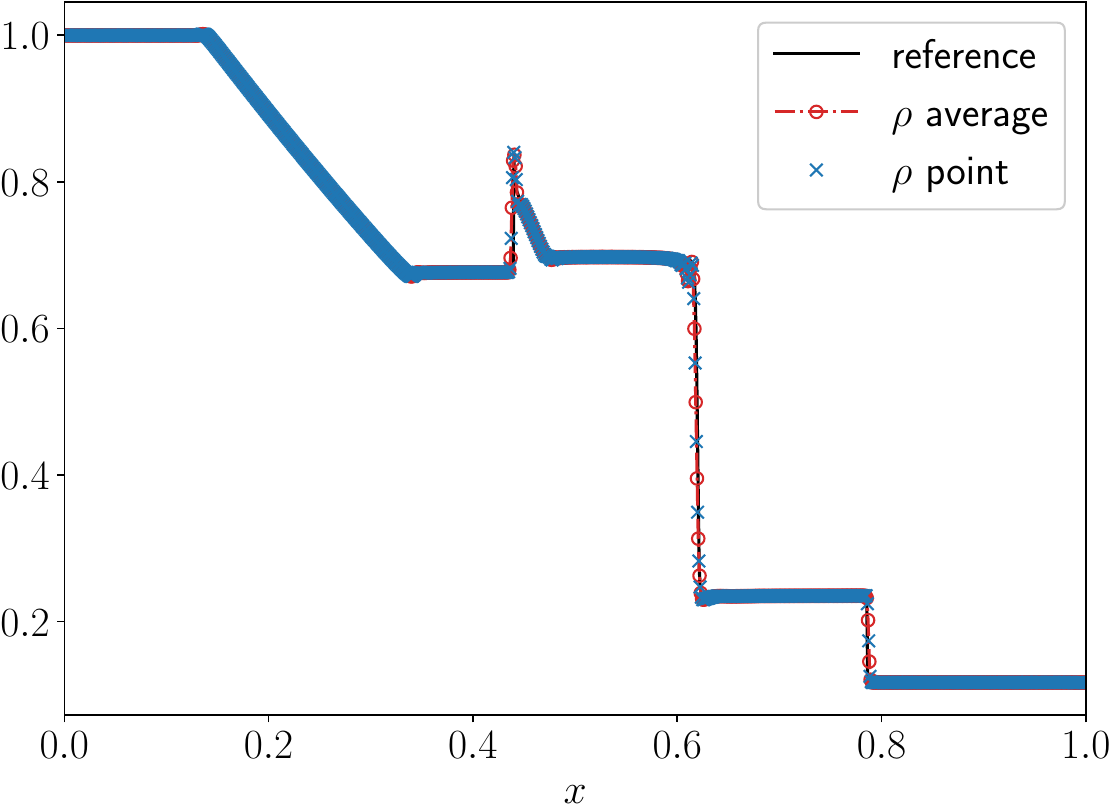}
		\end{subfigure}\hfill
		\begin{subfigure}[b]{0.32\textwidth}
			\centering
			\includegraphics[width=1.0\linewidth]{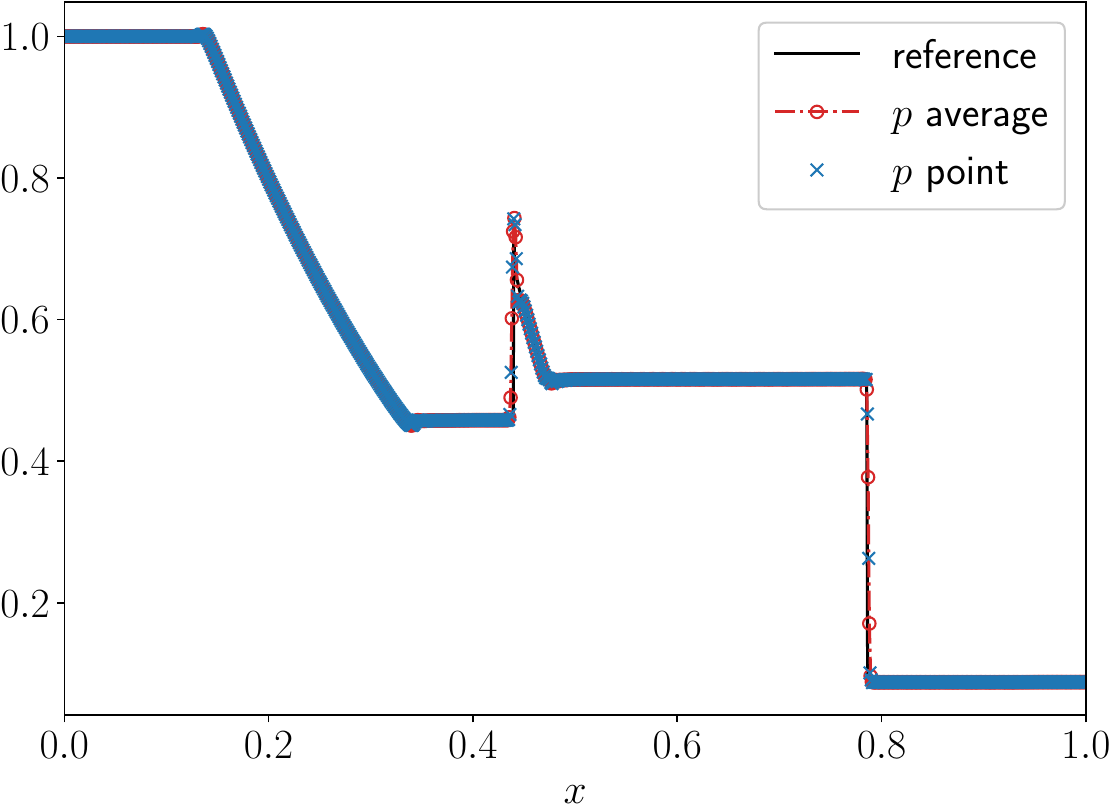}
		\end{subfigure}\hfill
		\begin{subfigure}[b]{0.32\textwidth}
			\centering
			\includegraphics[width=1.0\linewidth]{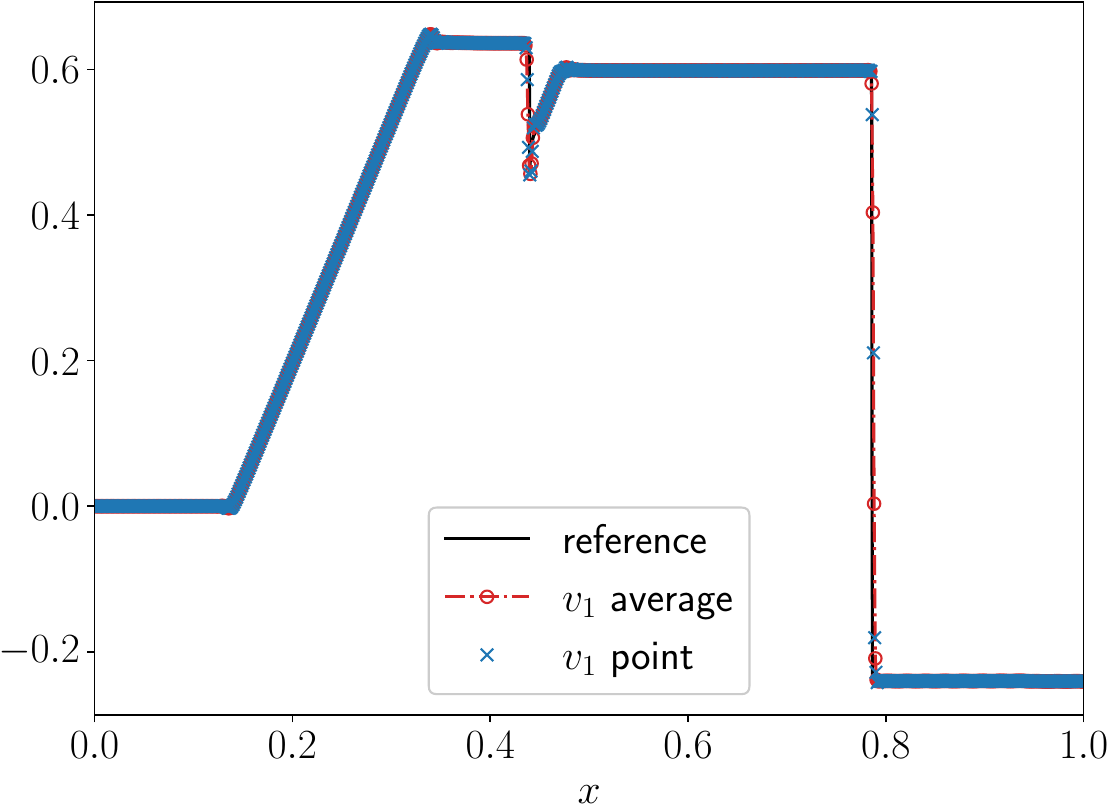}
		\end{subfigure}

		\medskip
		\begin{subfigure}[b]{0.32\textwidth}
			\centering
			\includegraphics[width=1.0\linewidth]{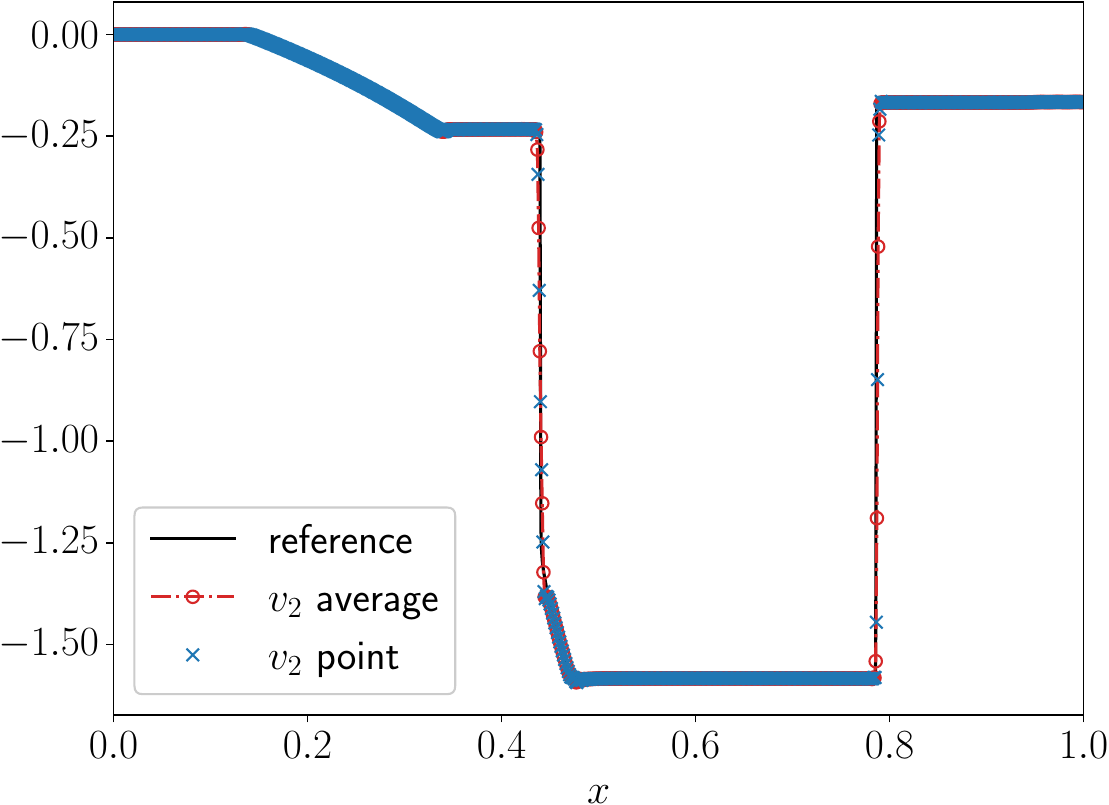}
		\end{subfigure}
		\begin{subfigure}[b]{0.32\textwidth}
			\centering
			\includegraphics[width=1.0\linewidth]{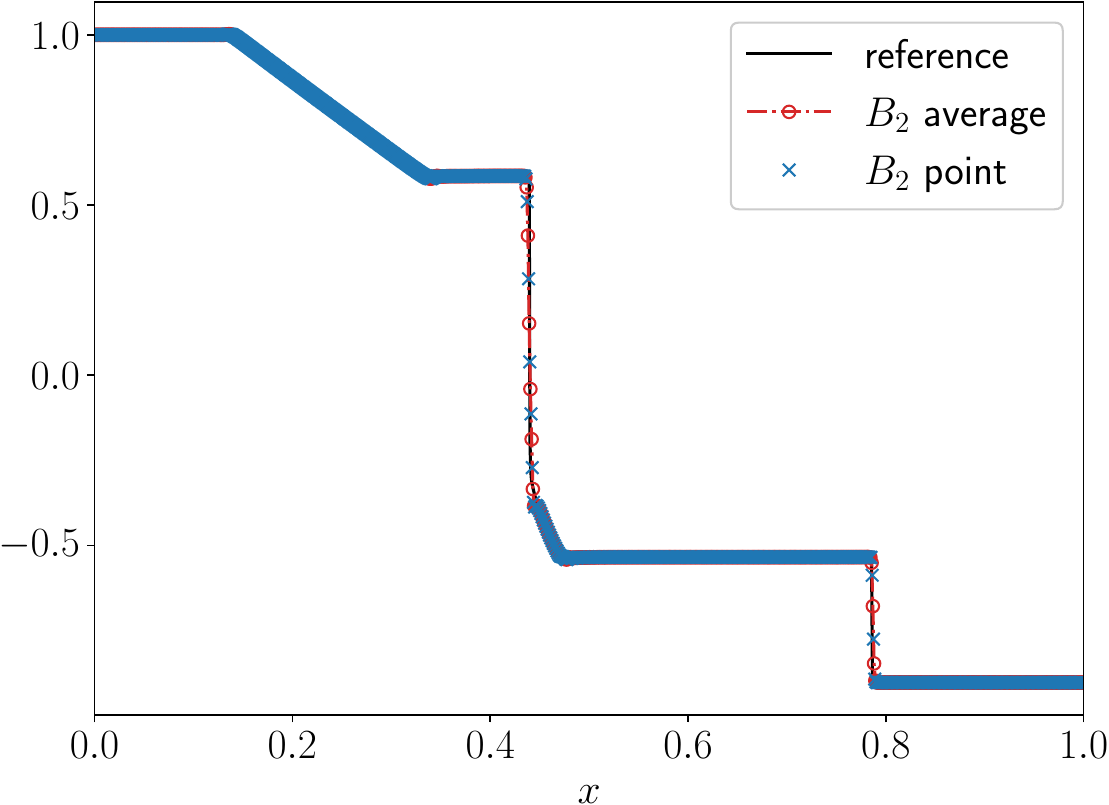}
		\end{subfigure}\hfill
		\caption{The first Riemann problem in Example \ref{ex:1d_rp},
			using $800$ cells and $\kappa=10$.}
		\label{fig:1d_rp1}
	\end{figure}
	
	\begin{figure}[htb!]
		\centering
		\begin{subfigure}[b]{0.325\textwidth}
			\centering
			\includegraphics[width=1.0\linewidth]{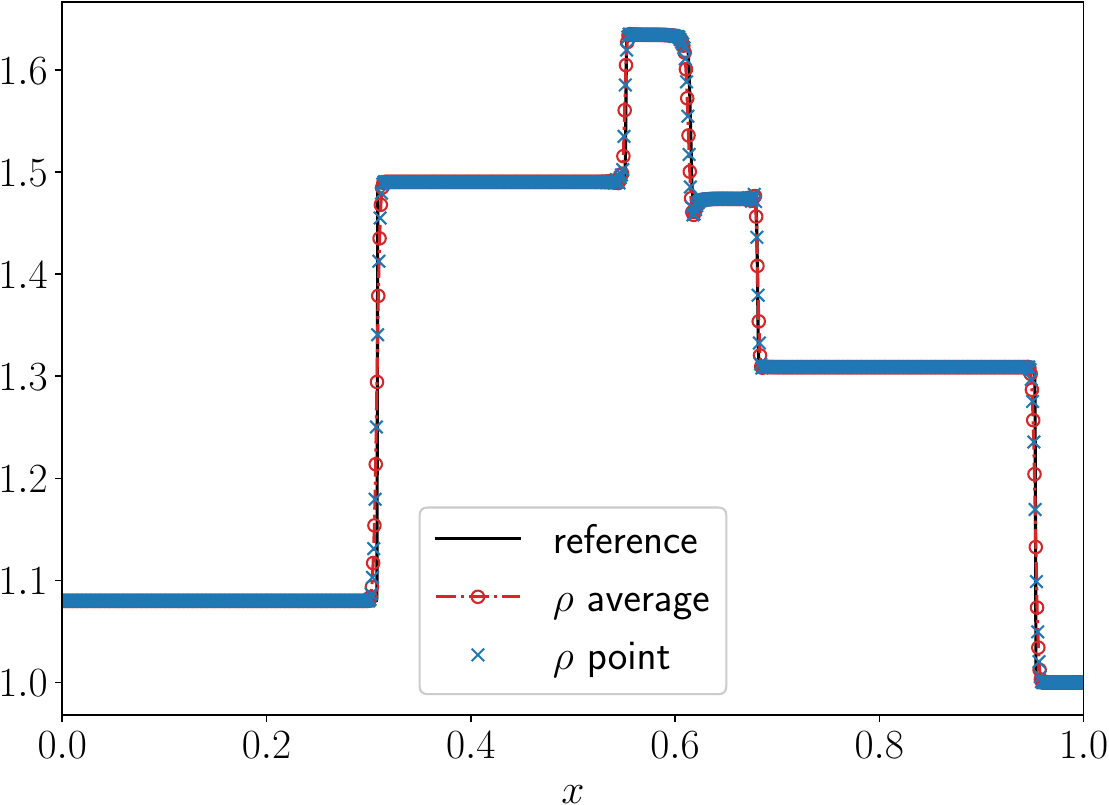}
		\end{subfigure}\hfill
		\begin{subfigure}[b]{0.325\textwidth}
			\centering
			\includegraphics[width=1.0\linewidth]{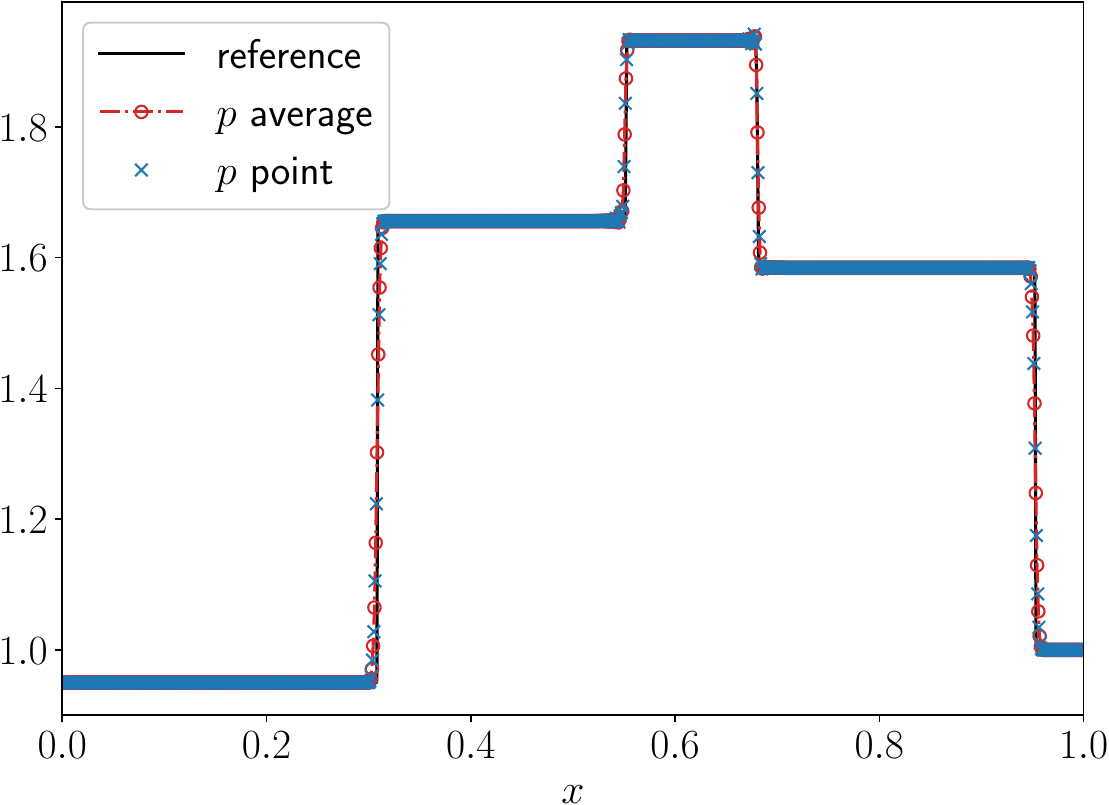}
		\end{subfigure}\hfill
		\begin{subfigure}[b]{0.325\textwidth}
			\centering
			\includegraphics[width=1.0\linewidth]{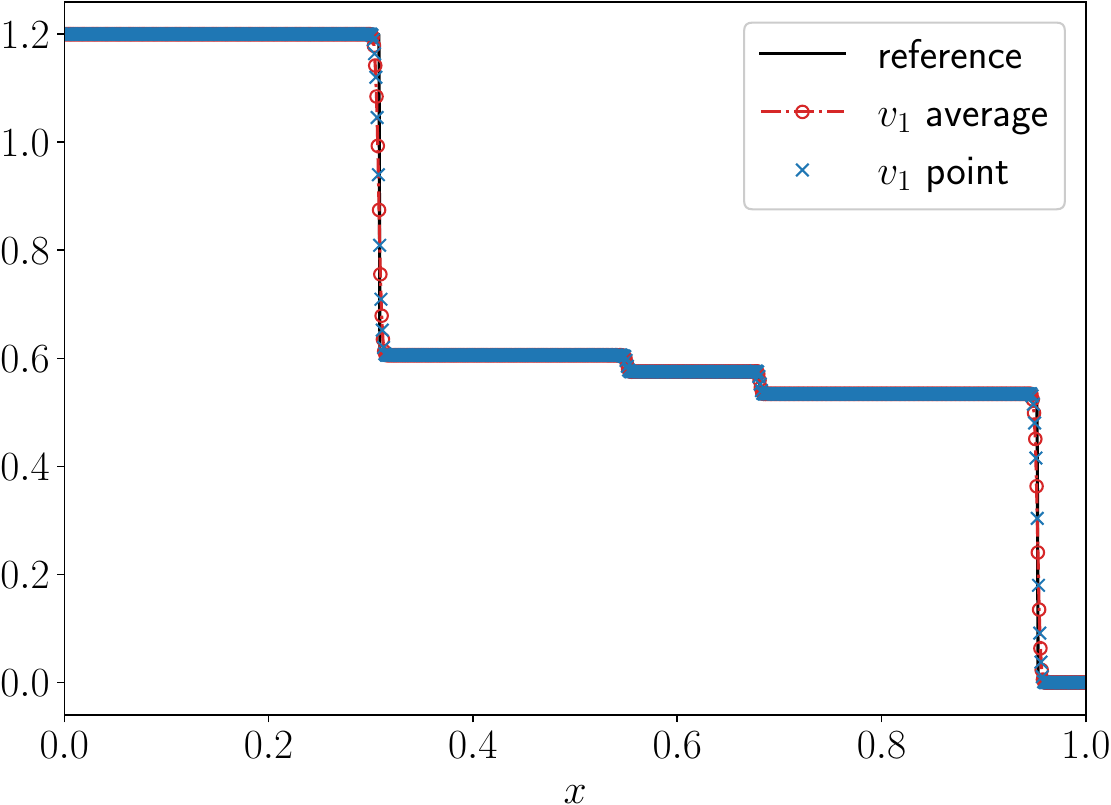}
		\end{subfigure}
		
		\medskip
		\begin{subfigure}[b]{0.325\textwidth}
			\centering
			\includegraphics[width=1.0\linewidth]{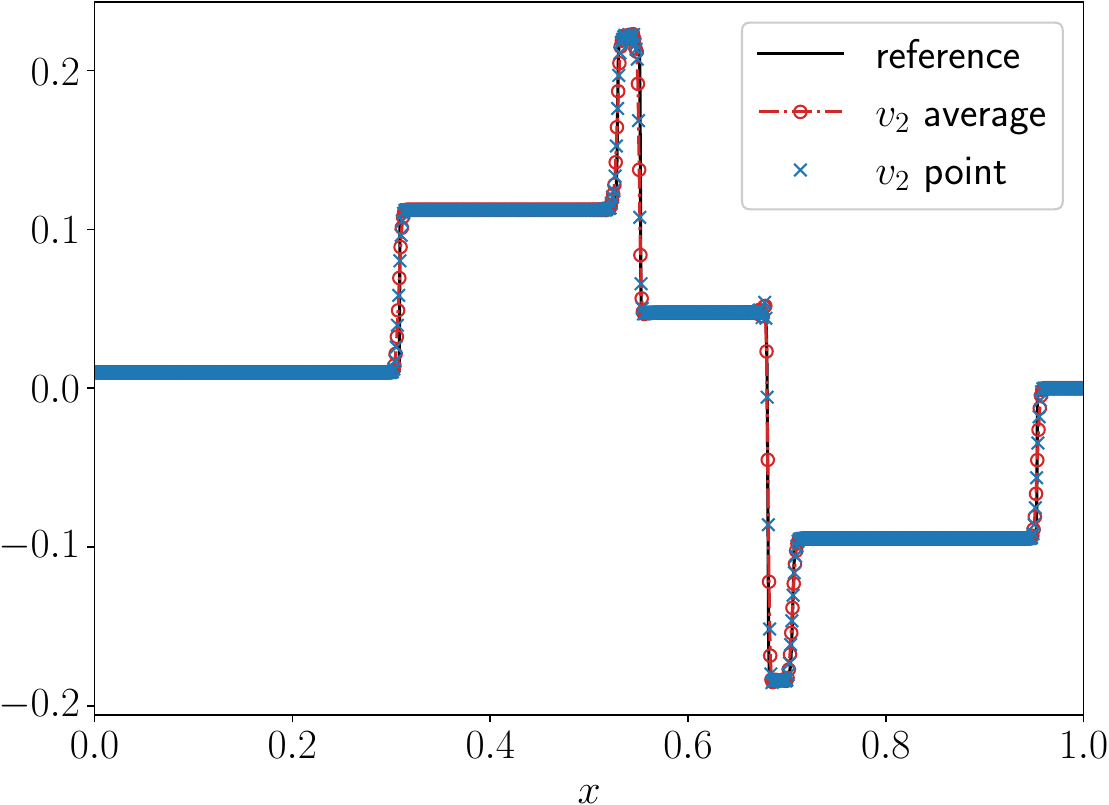}
		\end{subfigure}\hfill
		\begin{subfigure}[b]{0.325\textwidth}
			\centering
			\includegraphics[width=1.0\linewidth]{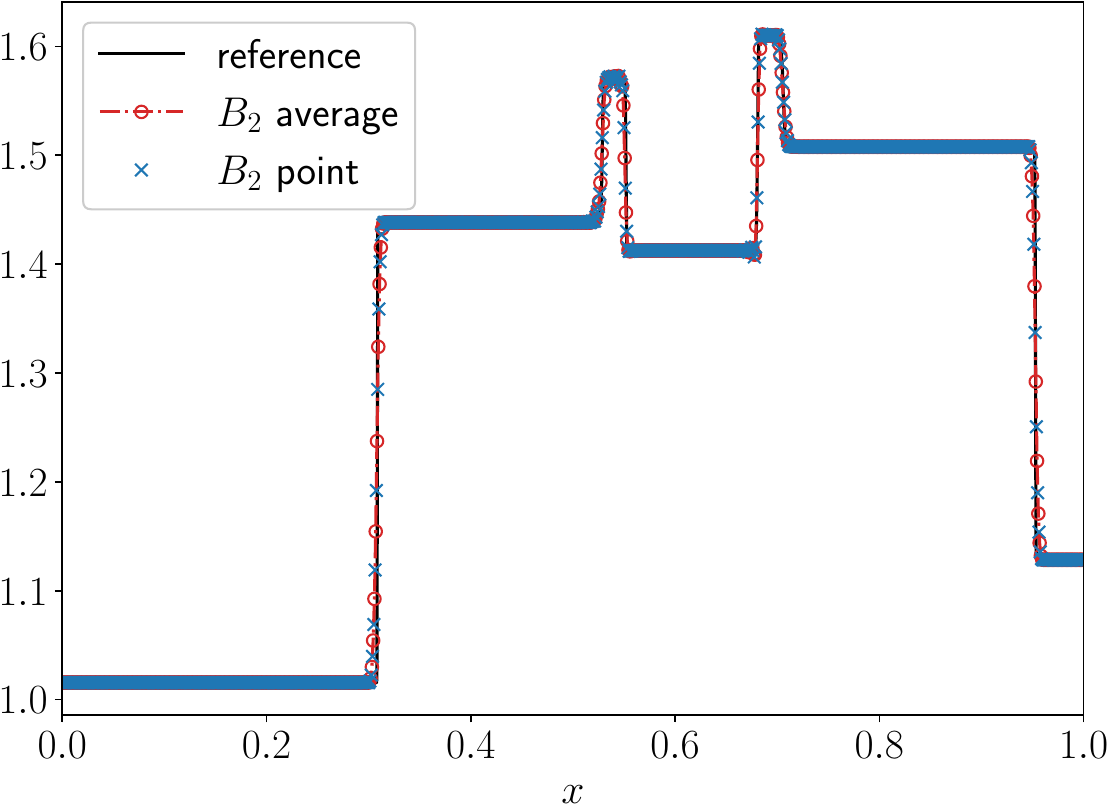}
		\end{subfigure}\hfill
		\begin{subfigure}[b]{0.325\textwidth}
			\centering
			\includegraphics[width=1.0\linewidth]{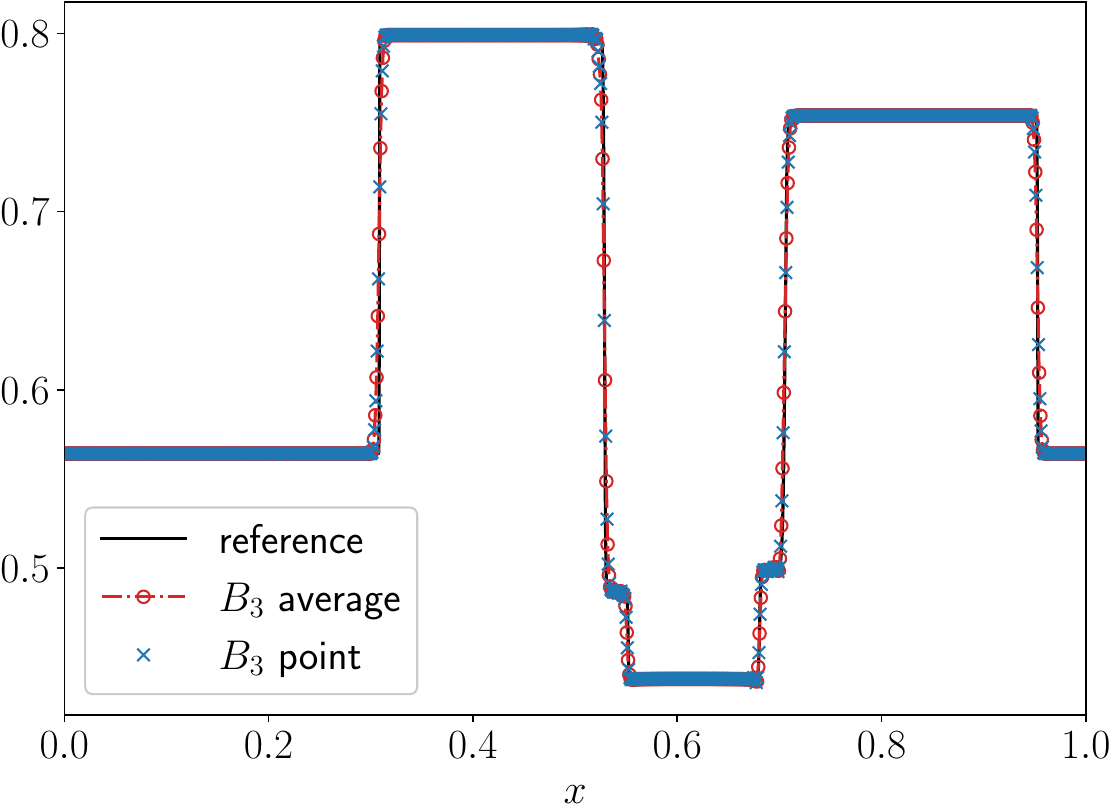}
		\end{subfigure}
		\caption{The second Riemann problem in Example \ref{ex:1d_rp},
		using $800$ cells and $\kappa=50$.}
		\label{fig:1d_rp2}
	\end{figure}
	
\end{example}

\begin{example}[1D Leblanc problem]\label{ex:1d_leblanc}
	To examine the PP property of our scheme, the test case in \cite{Liu_2025_Structure_JCP} is used.
	The computational domain is $[0,1]$ with the initial condition
	\begin{equation*}
		(\rho, \bv, \bB, p) =
		\begin{cases}
			(1, ~0, ~0, ~0, ~0.75, ~1, ~0, ~1), &\text{if}~x < 0.5, \\
			(0.125, ~0, ~0, ~0, ~0.75, ~-1, ~0, ~0.1), &\text{otherwise},
		\end{cases}
	\end{equation*}
	and the adiabatic index is $\gamma=1.4$.
	The final time is $T=1.5\times10^{-6}$.
	
	Figure \ref{fig:1d_leblanc} shows the density logarithm and magnetic pressure obtained with $2000$ cells and $\kappa=1000$,
	where the reference solution is obtained with $10^4$ cells.
	It is seen that the strong shock wave can be well captured without obvious oscillations.
	If the PP limitings are not activated, the simulation stops at the first time step due to negative pressure in point value, even if a small time step size $10^{-13}$ is used, which demonstrates the necessity of our PP limitings.
			
	\begin{figure}[htb!]
		\begin{subfigure}[b]{0.47\textwidth}
			\centering
			\includegraphics[width=1.0\linewidth]{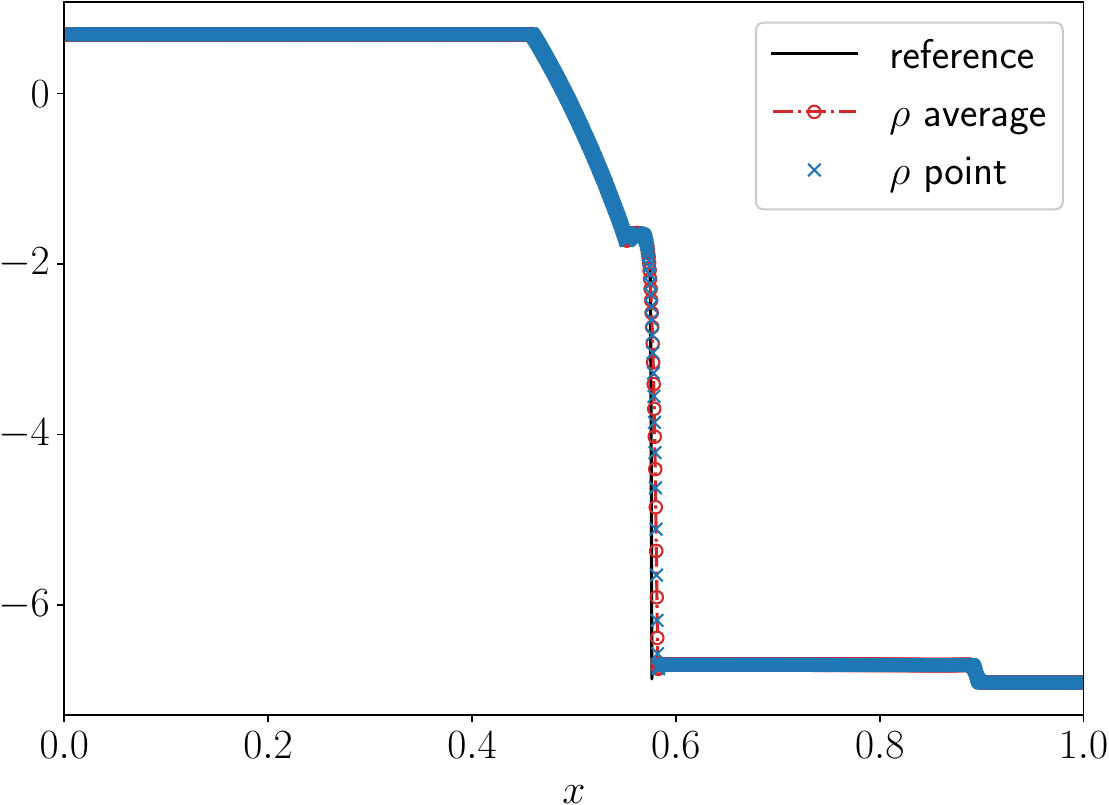}
		\end{subfigure}\hfill
		\begin{subfigure}[b]{0.493\textwidth}
			\centering
			\includegraphics[width=1.0\linewidth]{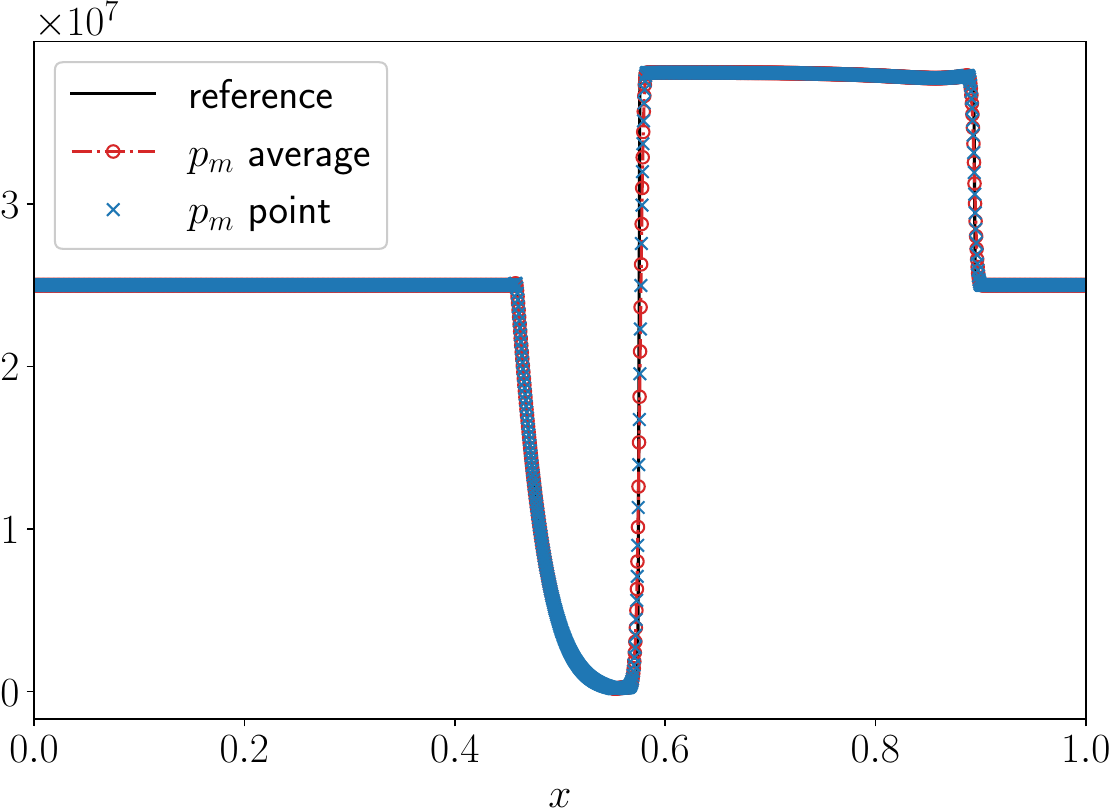}
		\end{subfigure}
		\caption{Example \ref{ex:1d_leblanc}.
		The density logarithm (left) and magnetic pressure (right).}
		\label{fig:1d_leblanc}
	\end{figure}
	
\end{example}

\begin{example}[2D accuracy tests]\label{ex:2d_accuracy}
	The adiabatic index is $\gamma=5/3$ for both cases.
	The first case is about a smooth sine wave propagating in the periodic domain $[0,1]\times[0,1]$ with the exact solution \cite{Wu_2018_provably_SJSC},
	\begin{equation*}
		(\rho, \bv, \bB, p) = (1 + 0.99\sin(2\pi(x+y-2t)), 1, 1, 0, 0.1, 0.1, 0, 1).
	\end{equation*}
	The errors in the $\ell_1$ norm at $T=0.1$ are shown in Figure \ref{fig:2d_accuracy}.
	The 3rd-order accuracy is obtained.
	
	In the second case, the MHD vortex problem \cite{Balsara_2004_Second_AJSS} is solved.
	The background flow $(\rho, \bv, \bB, p) = (1, 1, 1, 0, 0, 0, 0, 1)$ is initialized in the periodic domain $[-10, 10]\times[-10, 10]$.
	The following perturbation is added
	\begin{align*}
		(\delta v_1, \delta v_2) &= \xi\exp(0.5(1-r^2))(-y, x), \\
		(\delta B_1, \delta B_2) &= \mu\exp(0.5(1-r^2))(-y, x), \\
		\delta p &= 0.5(\mu^2(1-r^2) - \xi^2)\exp(1-r^2),
	\end{align*}
	where $r=\sqrt{x^2+y^2}$.
	The parameters are chosen as $\mu=5.389489439$, $\xi=\sqrt{2}\mu$ such that the lowest pressure at the vortex center is about $5.3\times 10^{-12}$ \cite{Christlieb_2015_Positivity_SJSC}.
	The PP limitings are necessary in this test, otherwise, the simulation stops due to negative pressure.
	The errors in the $\ell_1$ norm at $T=0.1$ are shown in Figure \ref{fig:2d_accuracy}, from which one observes almost 3rd-order accuracy.
	
	We also run this test with $\mu=1$ and without PP limitings by using upwind discretization based on the velocity direction for the source term in the point value update.
	The scheme is not stable, and generates negative pressure at around $T\approx 3.8885$ with $20\times20$ cells.
	When the central finite difference is used as in Section \ref{sec:pnt_evolution}, the 3rd-order accuracy can be obtained at $T=20$, which indicates that our point value update is stable.
	
	\begin{figure}[htb!]
		\begin{subfigure}[b]{0.48\textwidth}
			\centering
			\includegraphics[width=1.0\linewidth]{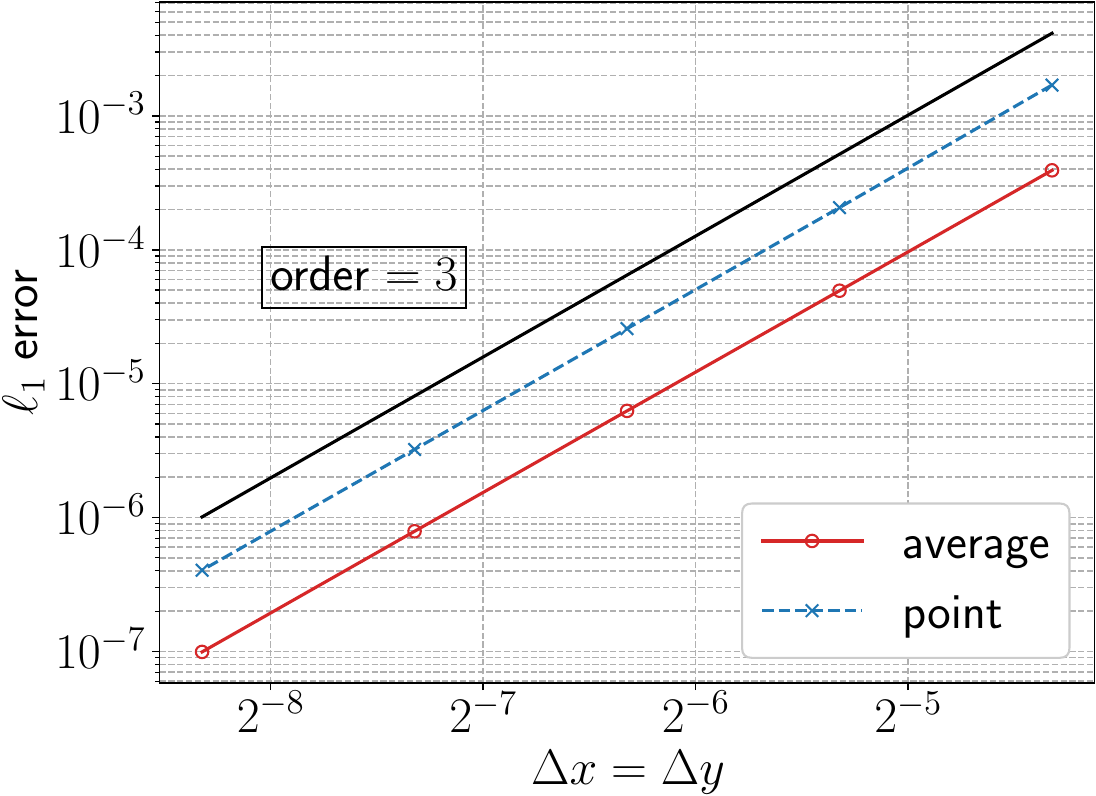}
		\end{subfigure}\hfill
		\begin{subfigure}[b]{0.48\textwidth}
			\centering
			\includegraphics[width=1.0\linewidth]{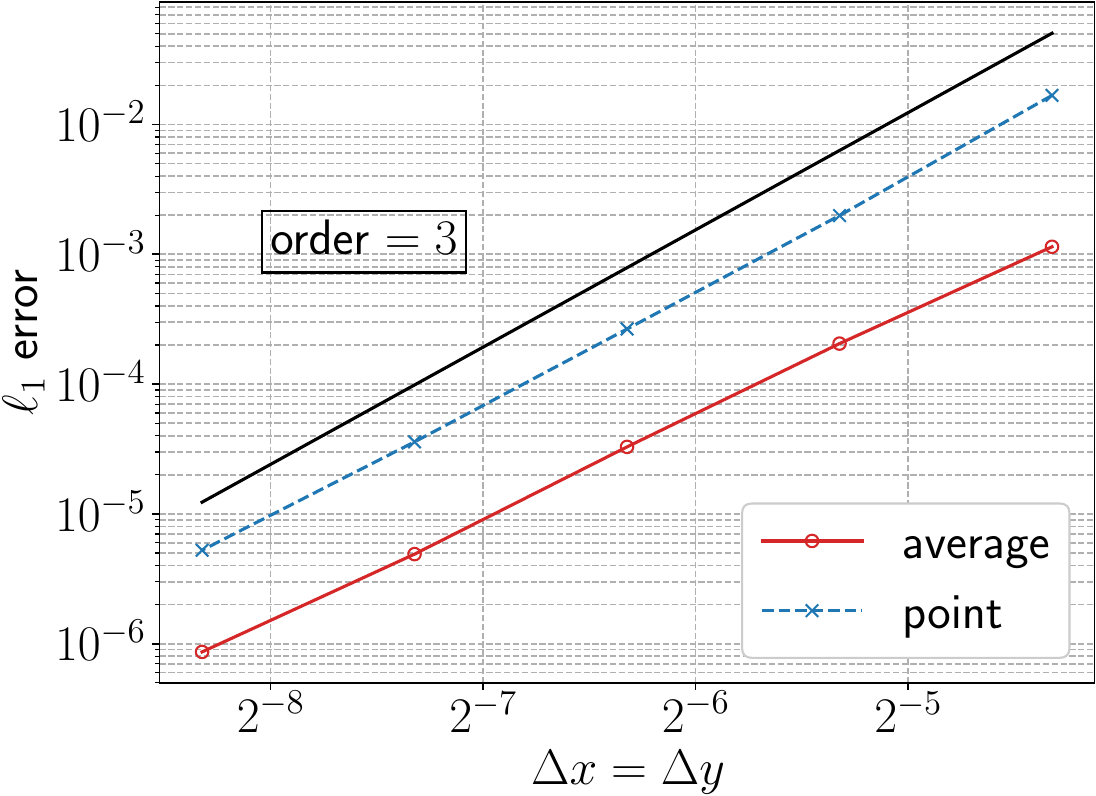}
		\end{subfigure}
		\caption{Example \ref{ex:2d_accuracy}.
			The errors and convergence rates of the smooth sine wave (left) and vortex (right).}
		\label{fig:2d_accuracy}
	\end{figure}
\end{example}

\begin{example}[Orszag-Tang problem]\label{ex:2d_orszag_tang}
	In this test, turbulent behavior will develop from smooth initial data  \cite{Orszag_1979_Small_JFM}.
	The domain is $[0,1]\times[0,1]$ with periodic boundary conditions, and the initial condition is	
	\begin{equation*}
		(\rho, \bv, \bB, p) = \left(\frac{25}{36\pi}, ~-\sin(2\pi y), ~\sin(2\pi x), ~0, ~-\frac{\sin(2\pi y)}{\sqrt{4\pi}}, ~\frac{\sin(4\pi x)}{\sqrt{4\pi}}, ~0, ~\frac{5}{12\pi}\right),
	\end{equation*}
	with the adiabatic index $\gamma=5/3$.
	
	The density plot at $T=0.5$ obtained by using $400\times400$ cells with $\kappa=1$ is shown in Figure \ref{fig:2d_orszag_tang} with the blending coefficients used in the shock sensor-based limiting.
	Our AF scheme can accurately capture the discontinuities and smooth structures, and the result is comparable to those in the literature.
    It is also seen that the shock sensor performs well.
	To examine the control of the divergence error, two kinds of discrete divergence are recorded.
	The first is computed based on the integration of the divergence in each cell, defined as
	\begin{equation*}
		\widetilde{(\divB)}_1(t) = \sum_{i,j} \sum_{l, m=1}^{3} \abs{(\divB)_{i,j}^{l, m}}\omega_l\omega_m\Delta x\Delta y,
	\end{equation*}
	where $(\divB)_{i,j}^{l, m}$ is obtained in the same way as \eqref{eq:semi_av_src}.
	The second is obtained by using the Gauss-Green formula, which approximates $\int_{\partial I_{ij}} \bB\cdot\bm{n}~ \dd s$ in each cell, i.e., the integration along edges, defined as
	\begin{align*}
		\widetilde{(\divB)}_2(t) = \sum_{i,j} \frac{1}{6}\Bigg\lvert &\big[(B_1)_{i+\frac12,j-\frac12} + 4(B_1)_{i+\frac12,j} + (B_1)_{i+\frac12,j+\frac12}\big] \\
		-& \big[(B_1)_{i-\frac12,j-\frac12} + 4(B_1)_{i-\frac12,j} + (B_1)_{i-\frac12,j+\frac12}\big] \\
		+& \big[(B_2)_{i-\frac12,j+\frac12} + 4(B_2)_{i,j+\frac12} + (B_2)_{i+\frac12,j+\frac12}\big] \\
		-& \big[(B_2)_{i-\frac12,j-\frac12} + 4(B_2)_{i,j-\frac12} + (B_2)_{i+\frac12,j-\frac12}\big] \Bigg\rvert \Delta x\Delta y.
	\end{align*}
	Figure \ref{fig:2d_orszag_tang_divB} shows their evolution in time.
	The discrete divergence is controlled in the sense that it increases very slowly over a long time, and almost arrives at a plateau.
	
	\begin{figure}[htb!]
		\begin{subfigure}[t]{0.325\textwidth}
			\centering
			\includegraphics[height=0.8\linewidth]{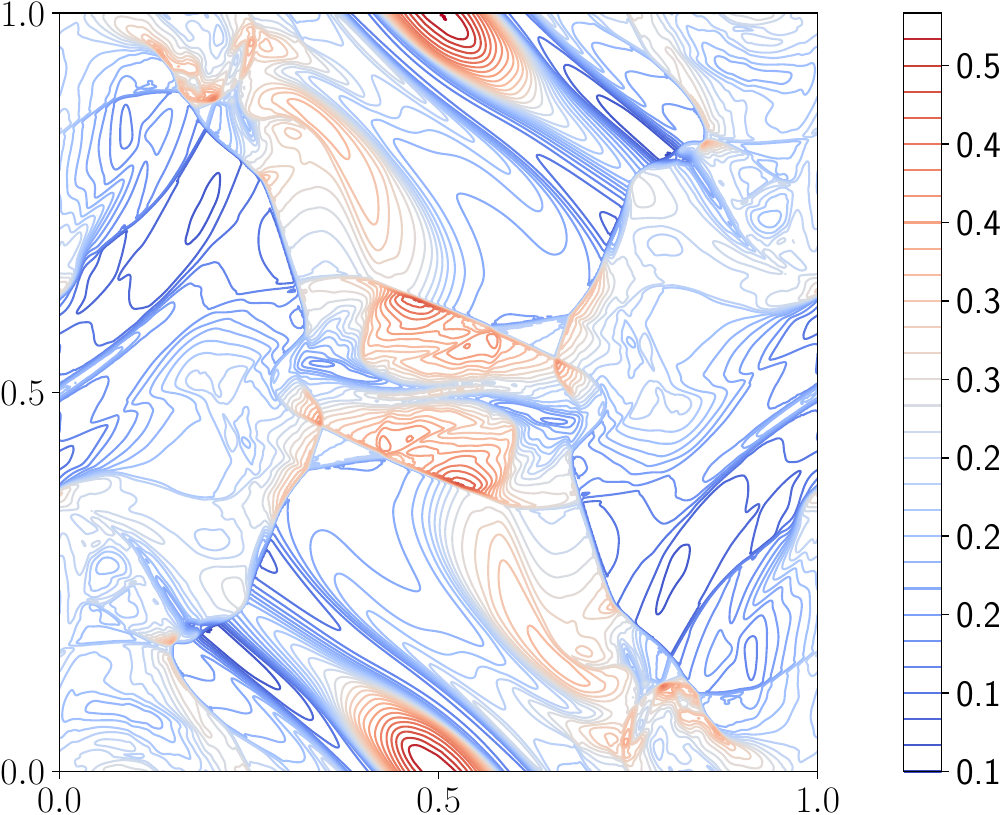}
		\end{subfigure}\hfill
		\begin{subfigure}[t]{0.325\textwidth}
			\centering
			\includegraphics[height=0.8\linewidth]{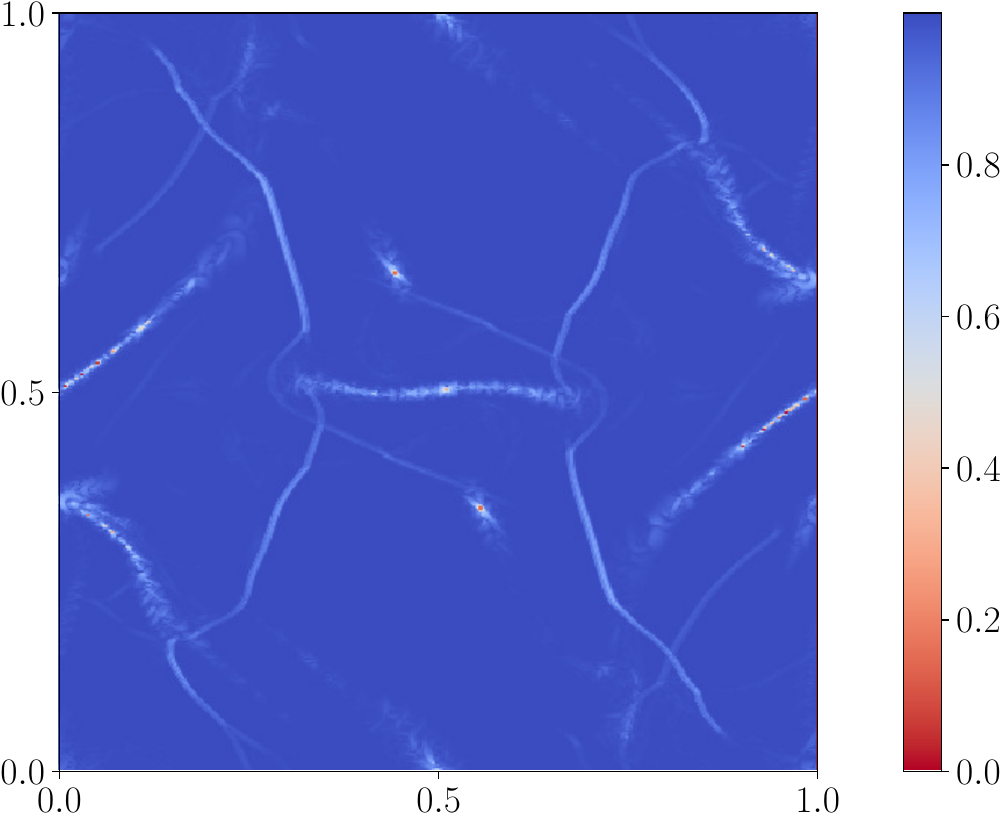}
		\end{subfigure}\hfill
		\begin{subfigure}[t]{0.325\textwidth}
			\centering
			\includegraphics[height=0.8\linewidth]{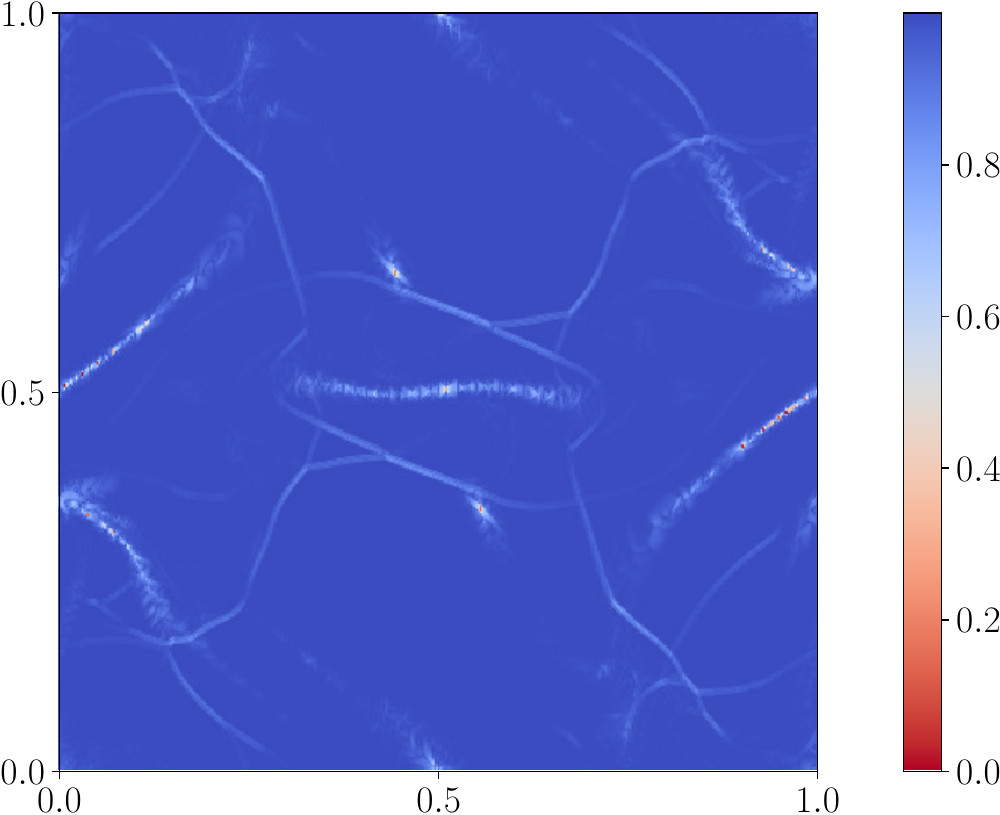}
		\end{subfigure}
		\caption{Example \ref{ex:2d_orszag_tang}.
			From left to right: $30$ equally spaced contour lines of the $\rho$ obtained by our PP AF scheme, the blending coefficients $\theta_{\xr,j}^s$, $\theta_{i,\yr}^s$ in the shock sensor.}
		\label{fig:2d_orszag_tang}
	\end{figure}
	
	\begin{figure}[htb!]
		\centering
		\includegraphics[width=0.6\linewidth]{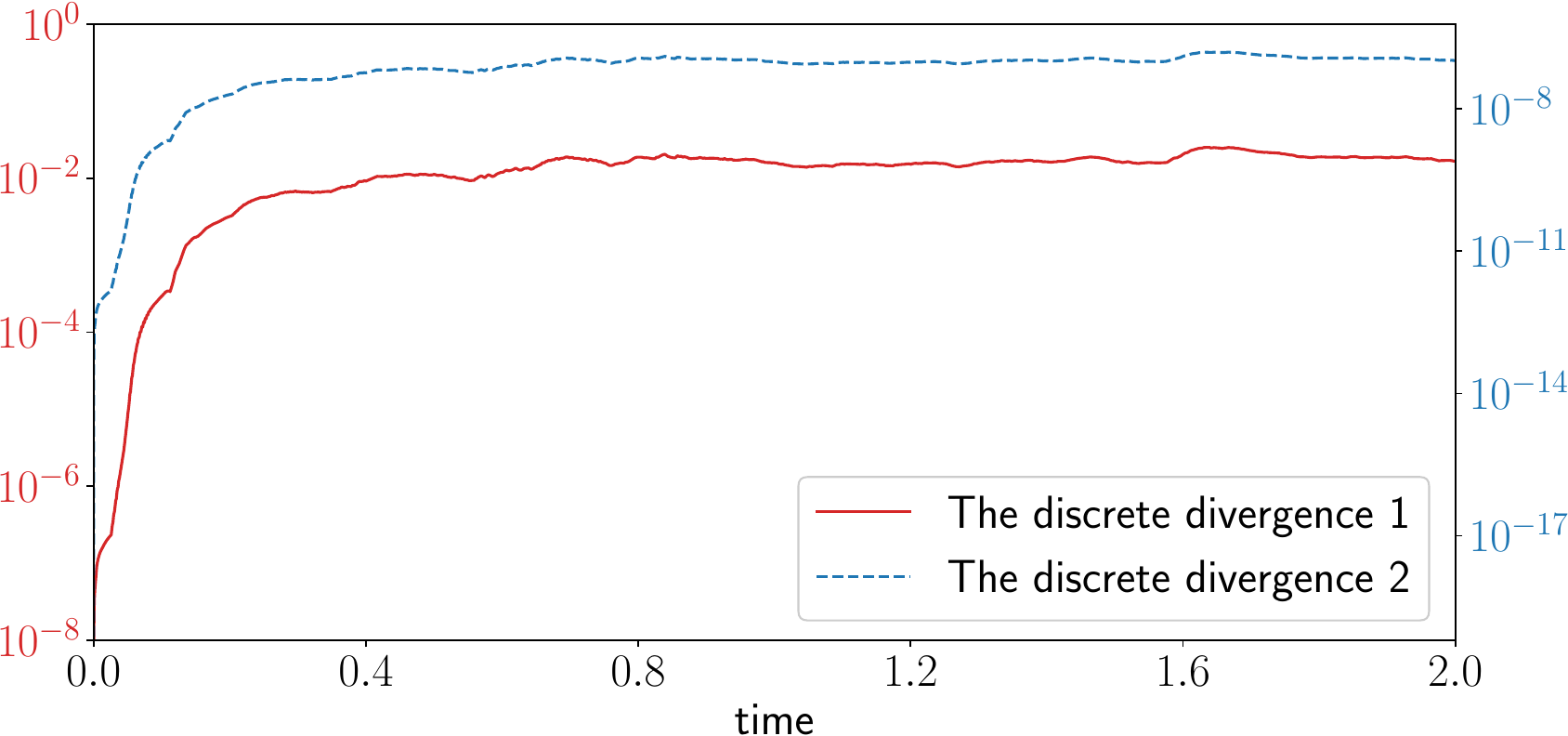}
		\caption{Example \ref{ex:2d_orszag_tang}.
			The evolution of the discrete divergence $\widetilde{(\divB)}_1(t)$ and $\widetilde{(\divB)}_2(t)$.}
		\label{fig:2d_orszag_tang_divB}
	\end{figure}
\end{example}

\begin{example}[Rotor problem]\label{ex:2d_rotor}
	This is the second rotor problem in \cite{Balsara_1999_staggered_JCP}, which describes a rotating dense fluid disk centered at a static background in the periodic domain $[0,1]\times[0,1]$.
	The magnetic field is initialized in the $x$-direction as $B_1=2.5/\sqrt{4\pi}$ and the pressure is $p=0.5$.
	The other initial data are
	\begin{equation*}
		(\rho, v_1, v_2) = \begin{cases}
			(10, ~-(y-0.5)/r_0, ~(x-0.5)/r_0), &\text{if}~ r < r_0, \\
			(1 + 9f, ~-f(y-0.5)/r, ~f(x-0.5)/r), &\text{if}~ r_0 < r < r_1, \\
			(1, ~0, ~0), &\text{if}~ r > r_1, \\
		\end{cases}
	\end{equation*}
	where $r=\sqrt{(x-0.5)^2+(y-0.5)^2}$, $r_0=0.1$, $r_1=0.115$,
	and $f=(r_1-r)/(r-r_0)$ is a tapper function.
	The adiabatic index is $\gamma=5/3$ and the test is solved until $T=0.295$.
	
	The numerical solutions obtained by using $400\times400$ cells and $\kappa=2$ are shown in Figure \ref{fig:2d_rotor},
	which are in good agreement with those in the literature.
	The shock sensor-based limiting is only used near the central rotor and the circular shock wave.
	Note that if the PP limitings are not activated, negative pressure appears at $T=7.6\times 10^{-2}$.
	Figure \ref{fig:2d_rotor_zoom} plots the enlarged view of the Mach number in the domain center.
	The left one with the Godunov-Powell source terms activated for both the cell average and point value preserves the circular rotation pattern well,
	while large distortions can be observed when the source terms are not used at the same time.
	This indicates that the divergence error is controlled in our AF scheme, as large divergence errors may cause distortion in the contour lines \cite{Balsara_1999_staggered_JCP,Li_2005_Locally_JSC,Toth_2000_$B0$_JCP}.
	
	\begin{figure}[htb!]
		\begin{subfigure}[t]{0.325\textwidth}
			\centering
			\includegraphics[height=0.8\linewidth]{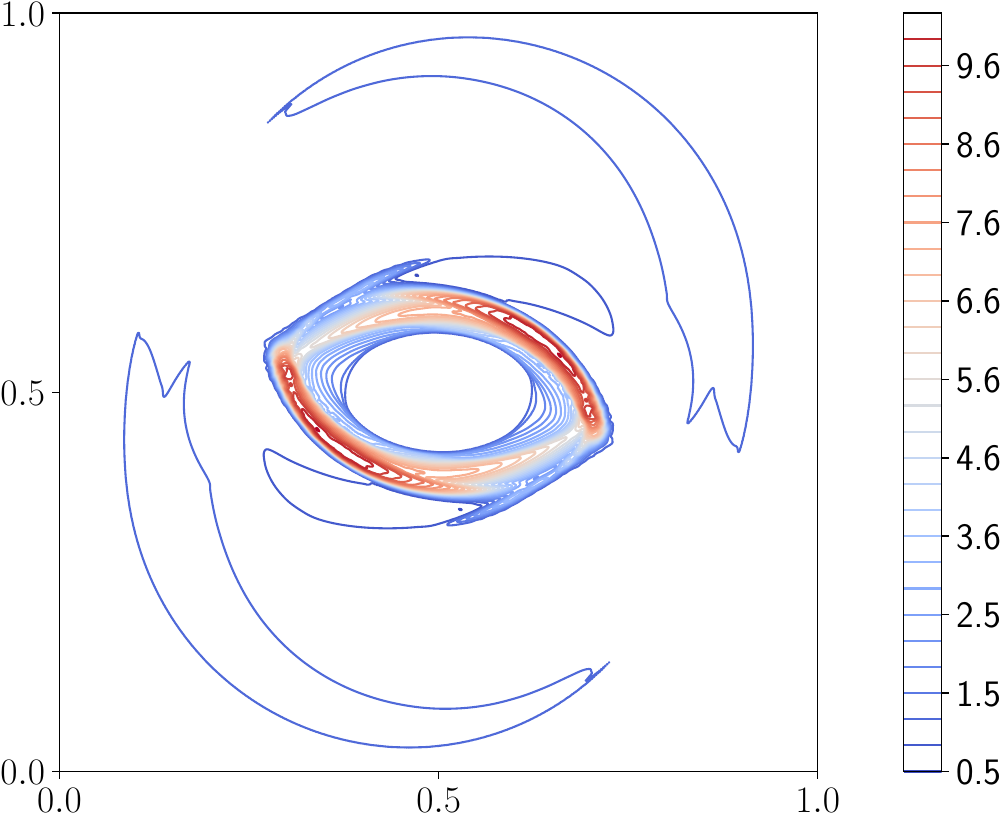}
			\caption{$\rho$}
		\end{subfigure}\hfill
		\begin{subfigure}[t]{0.325\textwidth}
			\centering
			\includegraphics[height=0.8\linewidth]{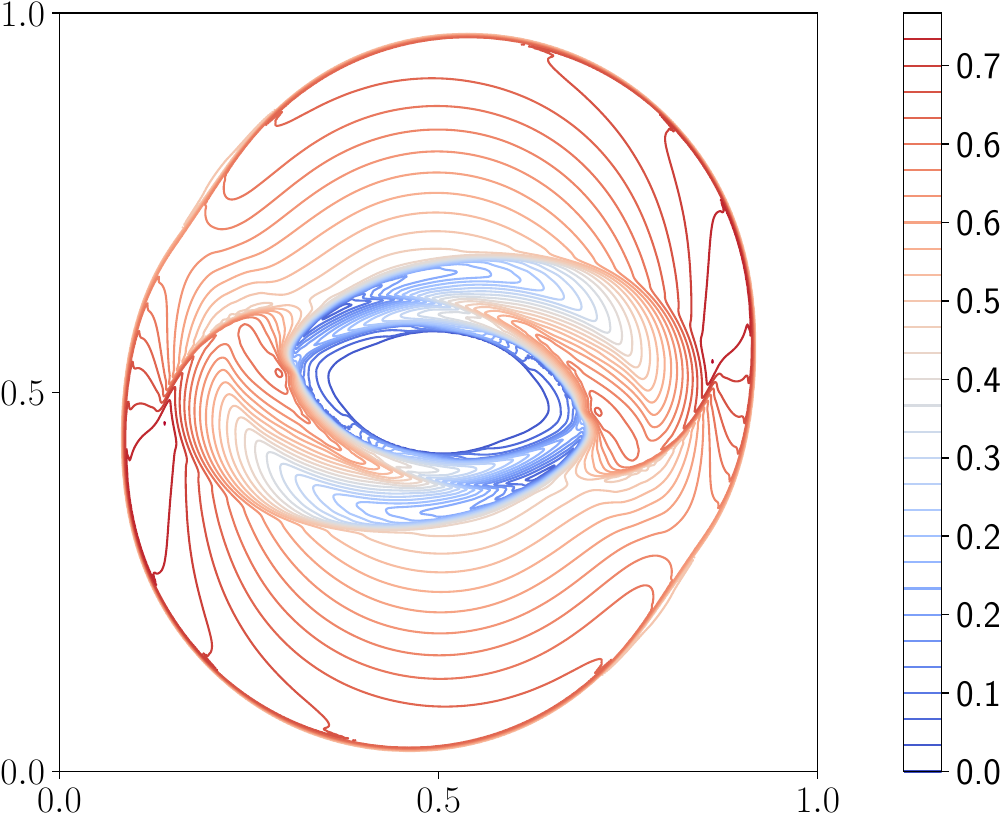}
			\caption{$p$}
		\end{subfigure}\hfill
		\begin{subfigure}[t]{0.325\textwidth}
			\centering
			\includegraphics[height=0.8\linewidth]{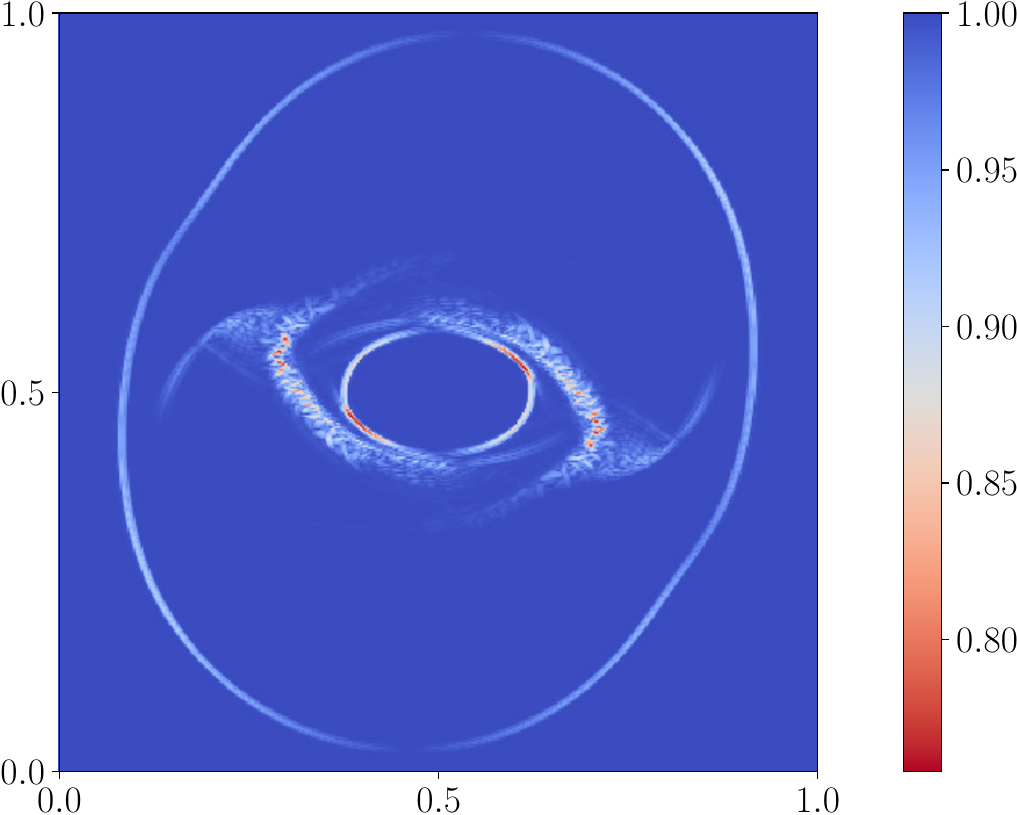}
			\caption{$\theta_{\xr,j}^s$}
		\end{subfigure}
		
		\begin{subfigure}[t]{0.325\textwidth}
			\centering
			\includegraphics[height=0.8\linewidth]{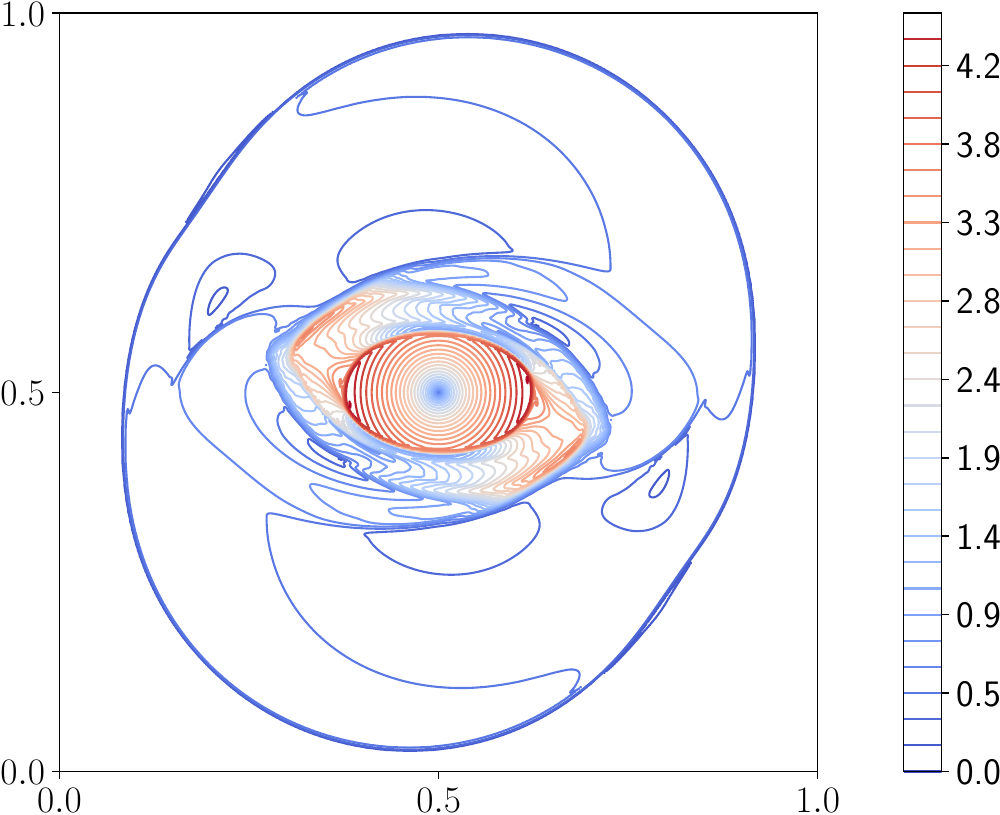}
			\caption{Mach number $\norm{\bv}/c$}
		\end{subfigure}\hfill
		\begin{subfigure}[t]{0.325\textwidth}
			\centering
			\includegraphics[height=0.8\linewidth]{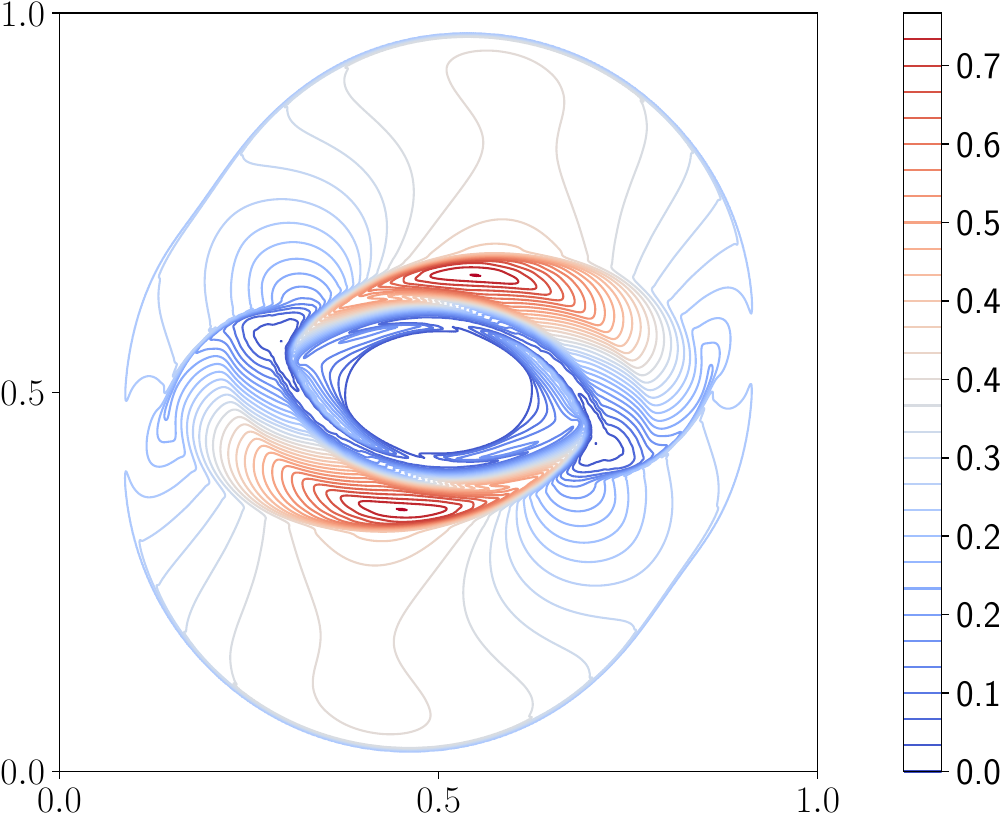}
			\caption{$p_m$}
		\end{subfigure}\hfill
		\begin{subfigure}[t]{0.325\textwidth}
			\centering
			\includegraphics[height=0.8\linewidth]{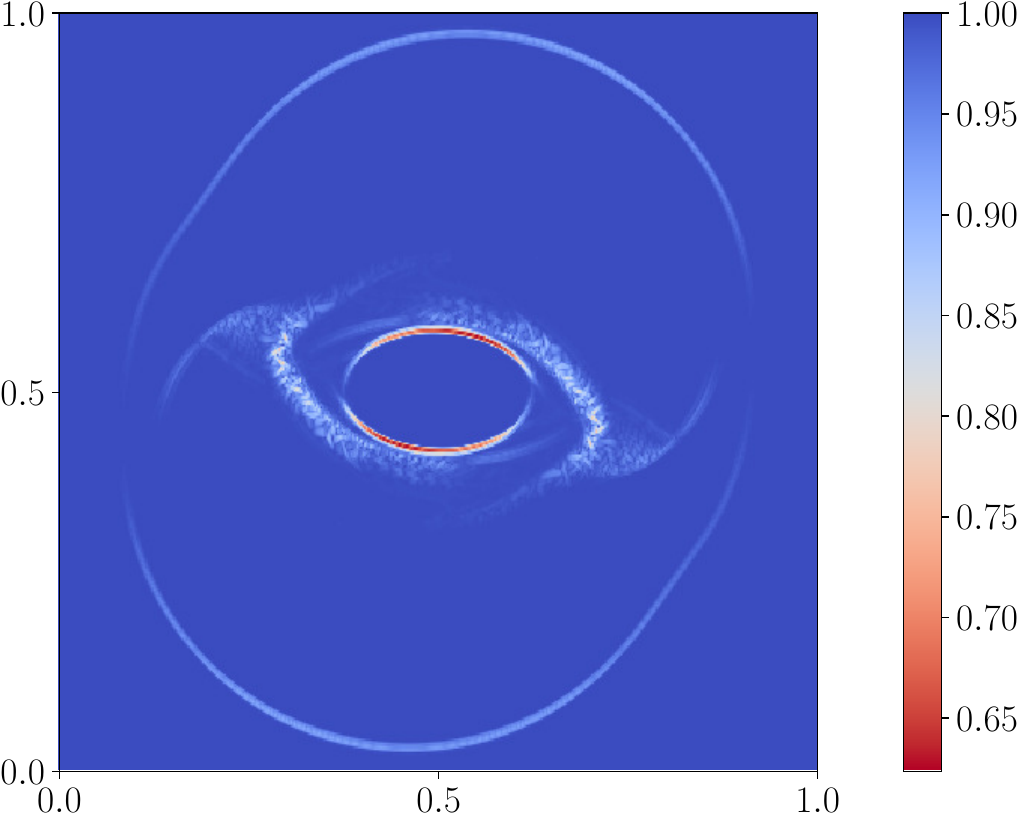}
			\caption{$\theta_{i,\yr}^s$}
		\end{subfigure}
		\caption{Example \ref{ex:2d_rotor}.
			$30$ equally spaced contour lines of the numerical solutions obtained by our PP AF scheme, and the blending coefficients in the shock sensor with $\kappa=2$.}
		\label{fig:2d_rotor}
	\end{figure}
	
	\begin{figure}[htb!]
		\begin{subfigure}[t]{0.24\textwidth}
			\centering
			\includegraphics[height=0.8\linewidth]{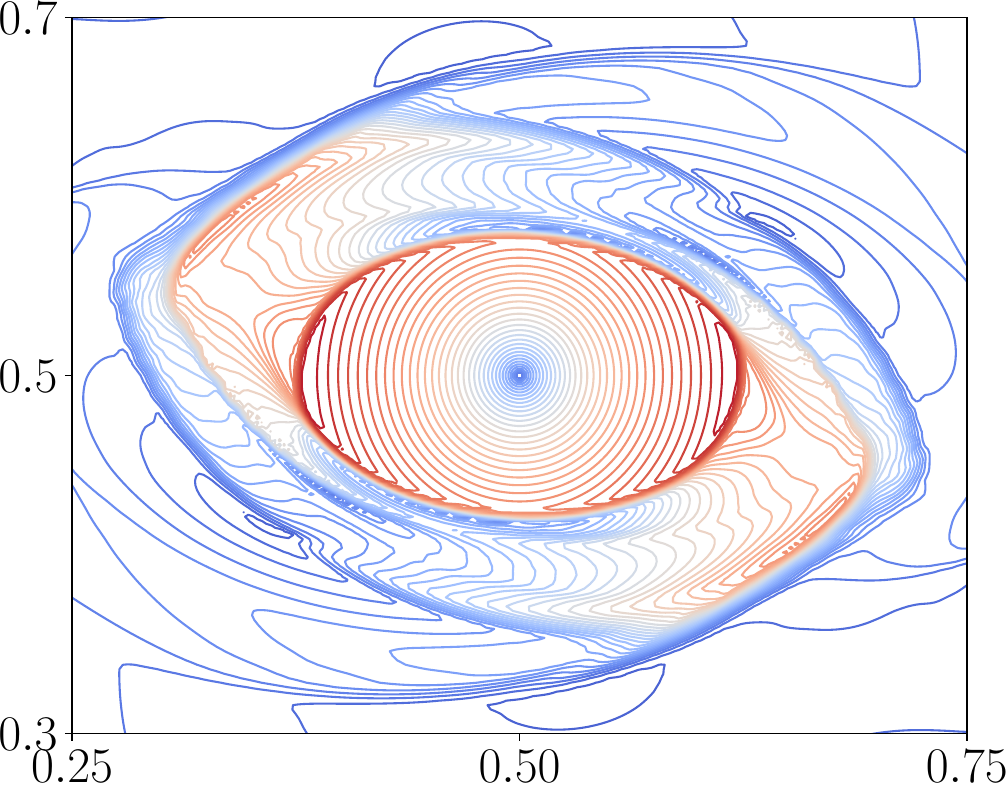}
		\end{subfigure}\hfill
		\begin{subfigure}[t]{0.24\textwidth}
			\centering
			\includegraphics[height=0.8\linewidth]{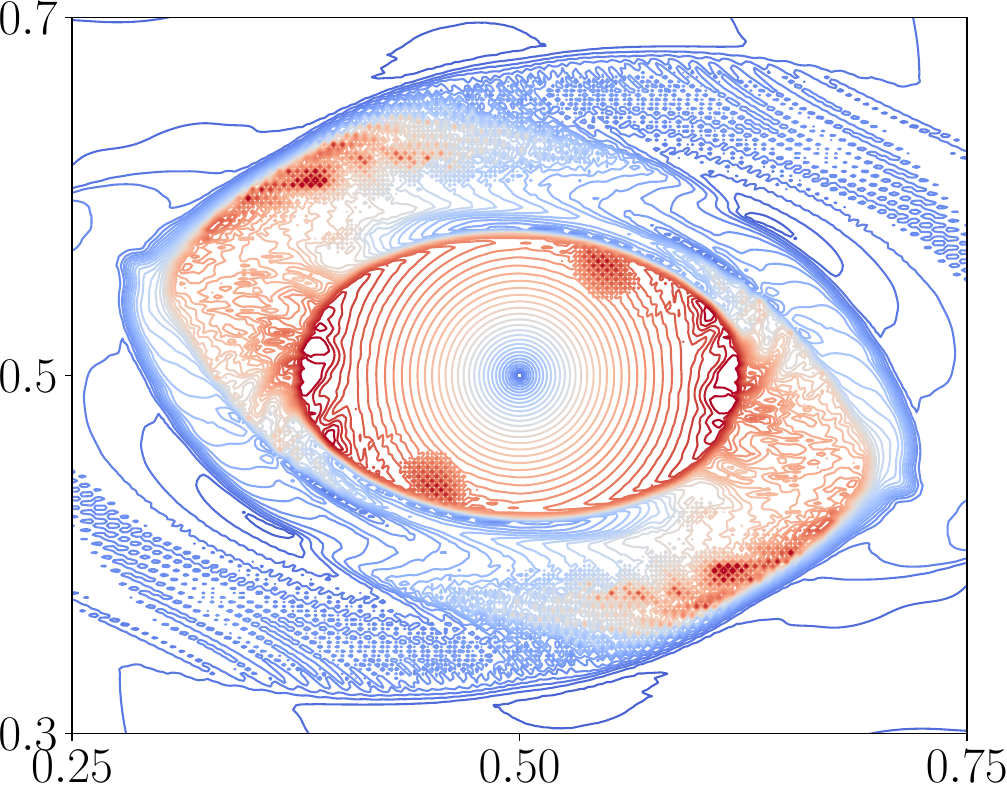}
		\end{subfigure}\hfill
		\begin{subfigure}[t]{0.24\textwidth}
			\centering
			\includegraphics[height=0.8\linewidth]{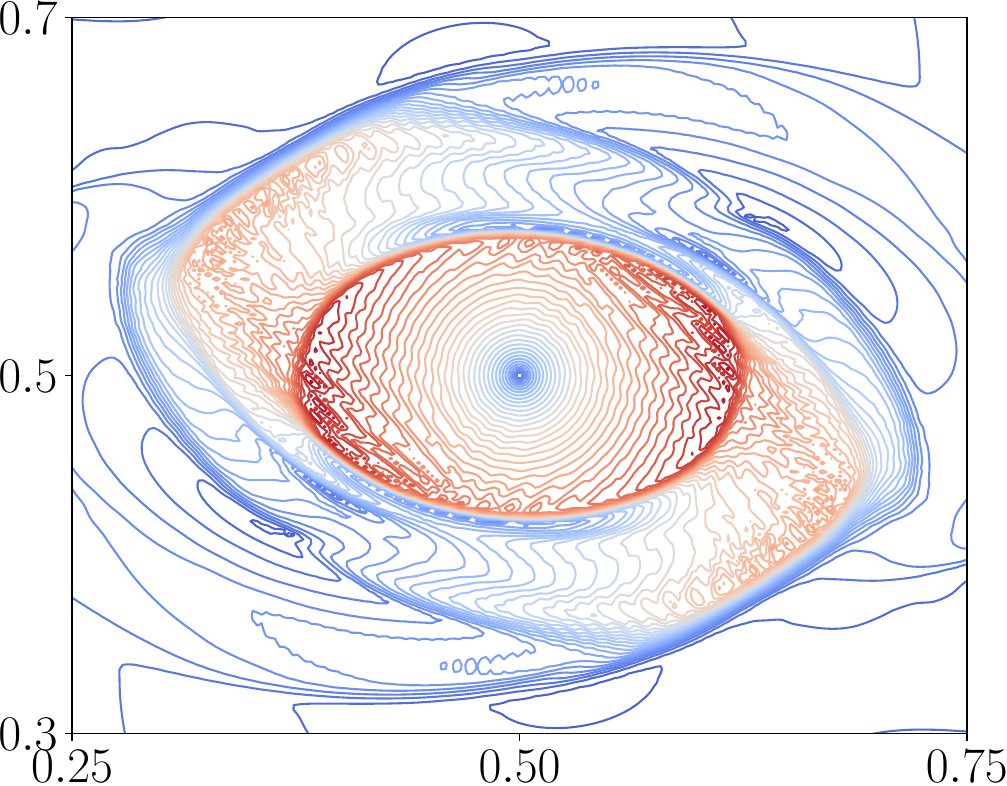}
		\end{subfigure}\hfill
		\begin{subfigure}[t]{0.24\textwidth}
			\centering
			\includegraphics[height=0.8\linewidth]{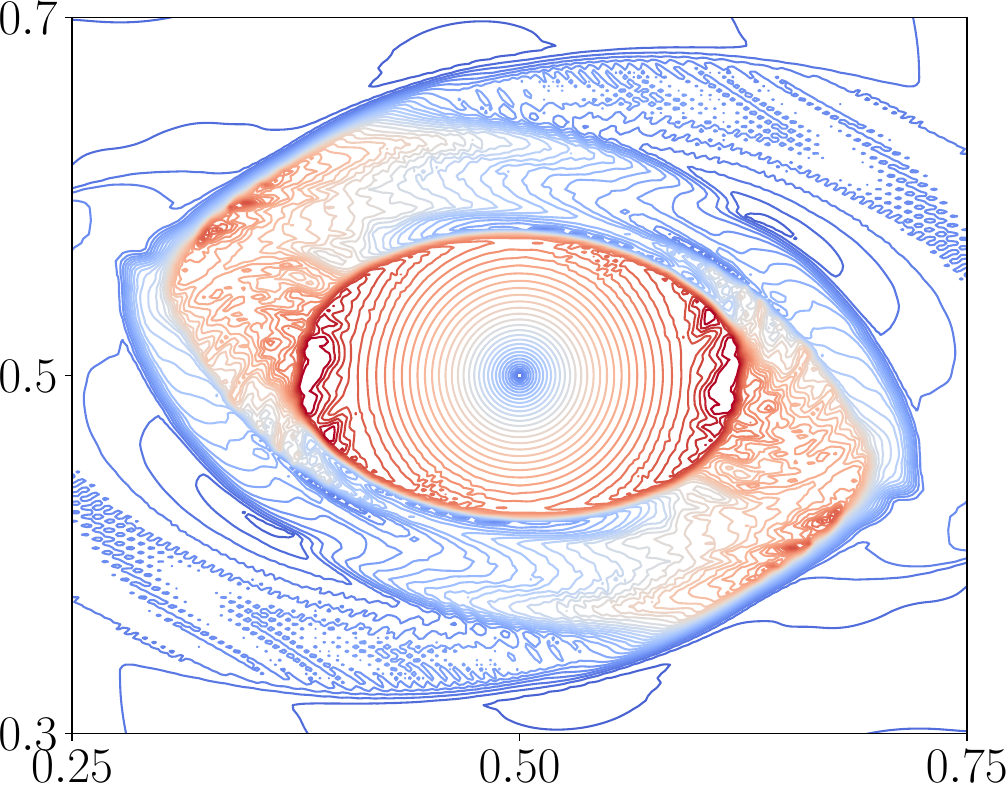}
		\end{subfigure}
		\caption{Example \ref{ex:2d_rotor}.
			$40$ equally spaced contour lines of the Mach number in the domain $[0.25,0.75]\times[0.3,0.7]$.
			From left to right: the source term is activated for both the cell average and point value, only for point value, only for cell average, neither.}
		\label{fig:2d_rotor_zoom}
	\end{figure}
		
\end{example}

\begin{example}[Blast wave]\label{ex:2d_blast}
	This is a test problem with a strongly magnetized medium with low plasma.
	Following the setup in \cite{Stone_2008_Athena_AJSS},
	the computational domain is $[-0.5,0.5]\times[-0.5,0.5]$ with outflow boundary conditions, and the initial condition is	
	\begin{equation*}
		(\rho, \bv, \bB, p) = (1, ~0, ~0, ~0, ~1/\sqrt{2}, ~1/\sqrt{2}, ~0, ~0.1),
	\end{equation*}
	except for a larger pressure $p=10$ in the central circular part $\sqrt{x^2+y^2}<0.1$.
	The adiabatic index is $\gamma=5/3$, and the problem is solved until $T=0.2$.
	
	The results obtained by using our PP AF scheme with $400\times400$ cells and $\kappa=1$ are shown in Figure \ref{fig:2d_blast}.
	The flow structures, including the outward-going circular blast wave, are captured with high resolution, which agree well with those in \cite{Stone_2008_Athena_AJSS,Liu_2025_Structure_JCP}.
	The simulation stops due to negative pressure if the PP limitings are not activated.
	One can also observe that the shock sensor-based limiting is activated locally.
	
	\begin{figure}[htb!]
		\begin{subfigure}[t]{0.325\textwidth}
			\centering
			\includegraphics[height=0.8\linewidth]{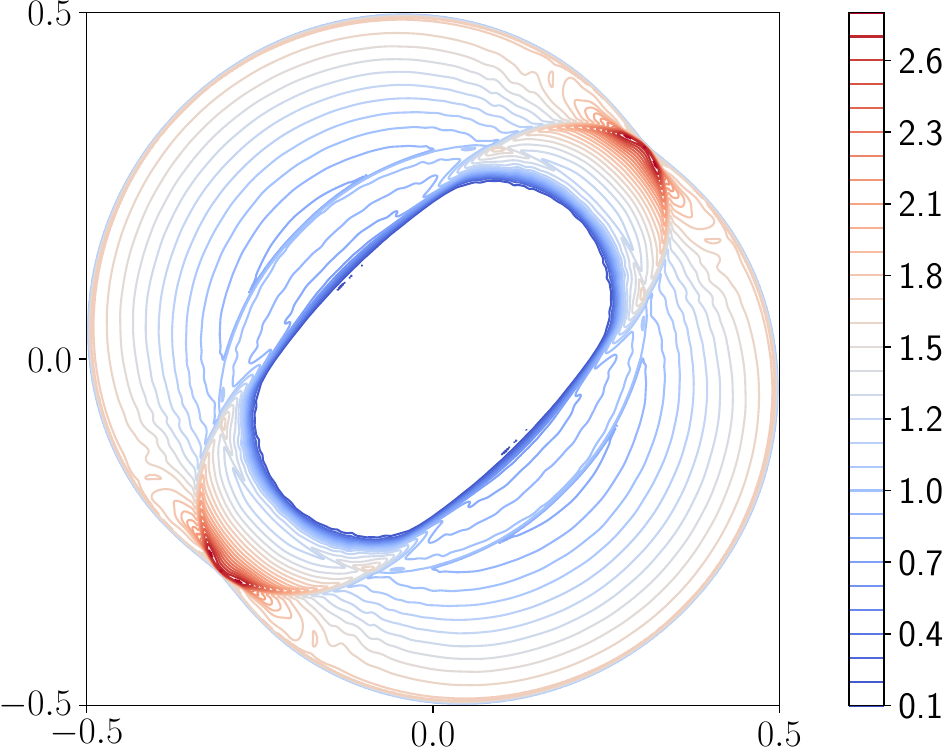}
			\caption{$\rho$}
		\end{subfigure}\hfill
		\begin{subfigure}[t]{0.325\textwidth}
			\centering
			\includegraphics[height=0.8\linewidth]{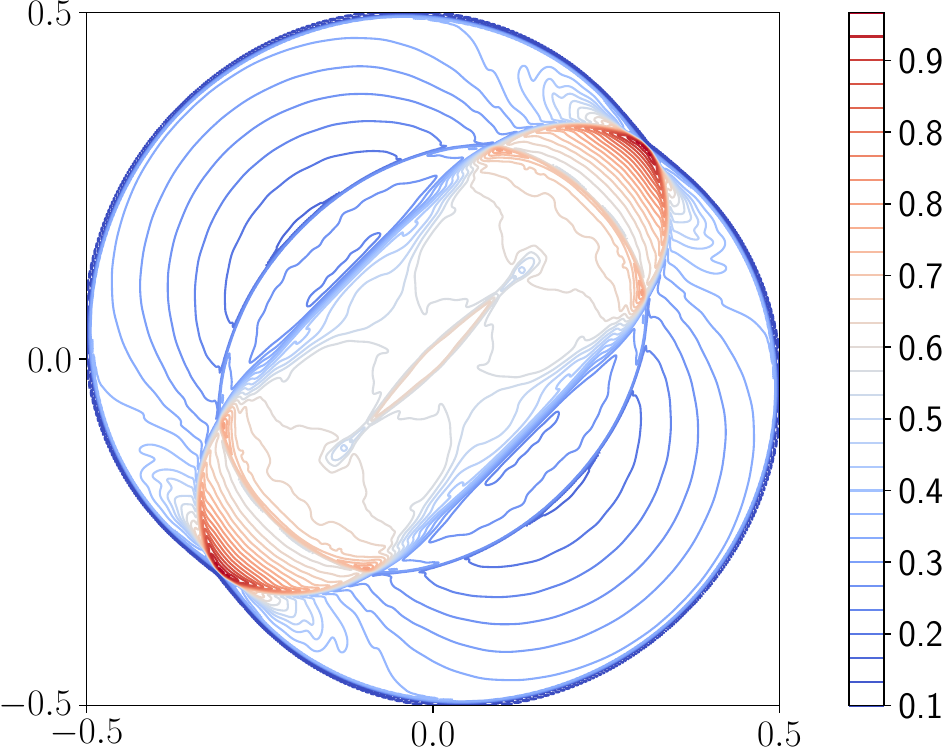}
			\caption{$p$}
		\end{subfigure}\hfill
		\begin{subfigure}[t]{0.325\textwidth}
			\centering
			\includegraphics[height=0.8\linewidth]{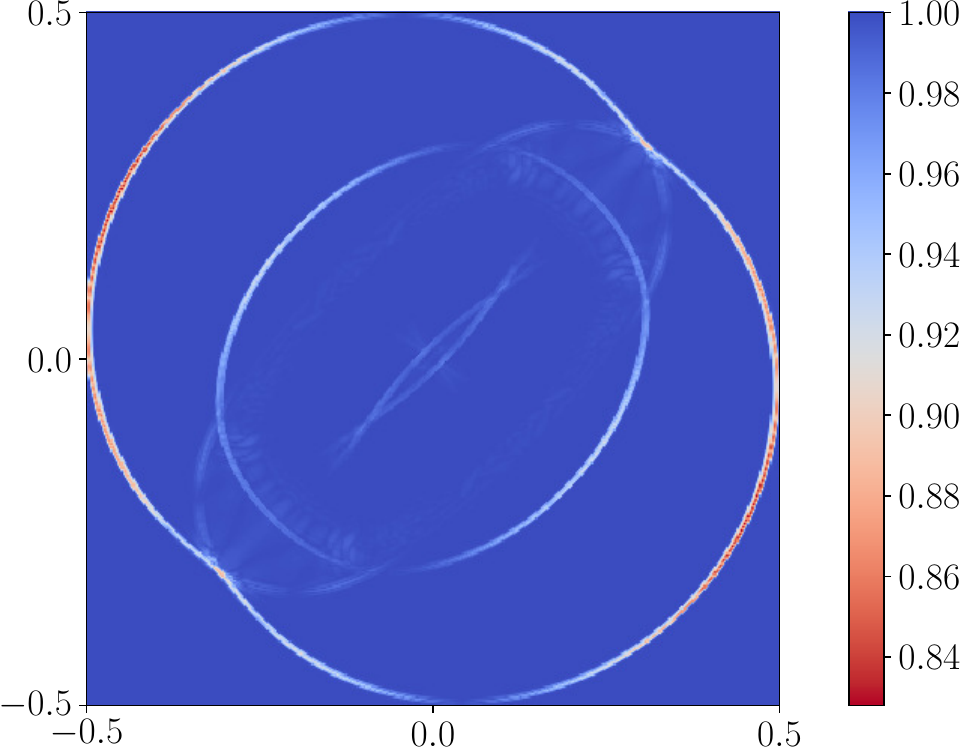}
			\caption{$\theta_{i+\frac12,j}^s$}
		\end{subfigure}
		
		\begin{subfigure}[t]{0.325\textwidth}
			\centering
			\includegraphics[height=0.8\linewidth]{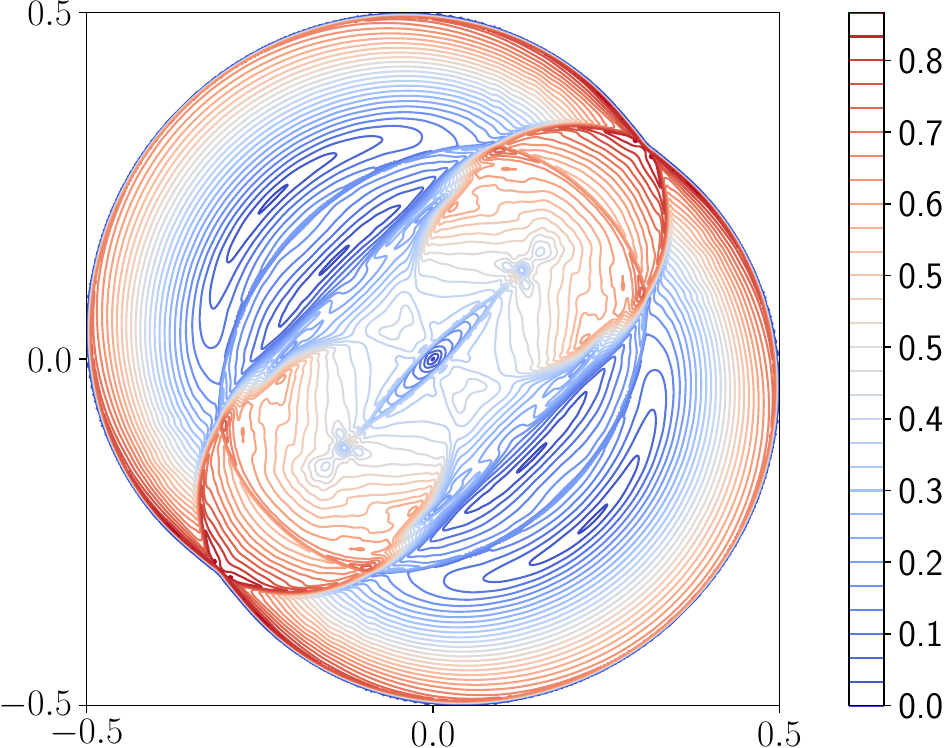}
			\caption{$\norm{\bv}$}
		\end{subfigure}\hfill
		\begin{subfigure}[t]{0.325\textwidth}
			\centering
			\includegraphics[height=0.8\linewidth]{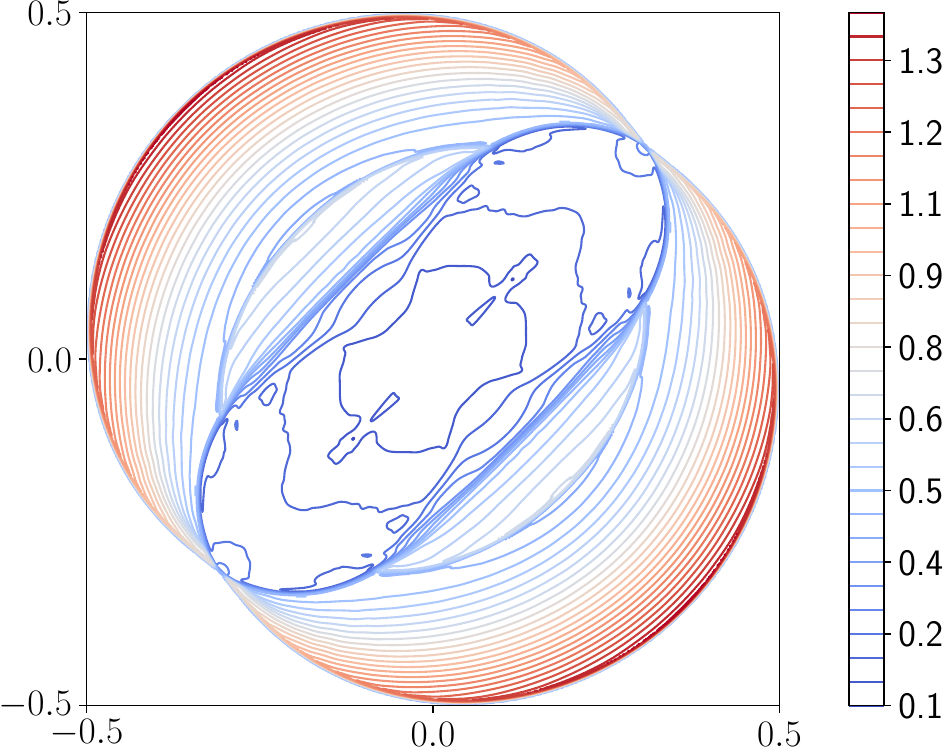}
			\caption{$p_m$}
		\end{subfigure}\hfill
		\begin{subfigure}[t]{0.325\textwidth}
			\centering
			\includegraphics[height=0.8\linewidth]{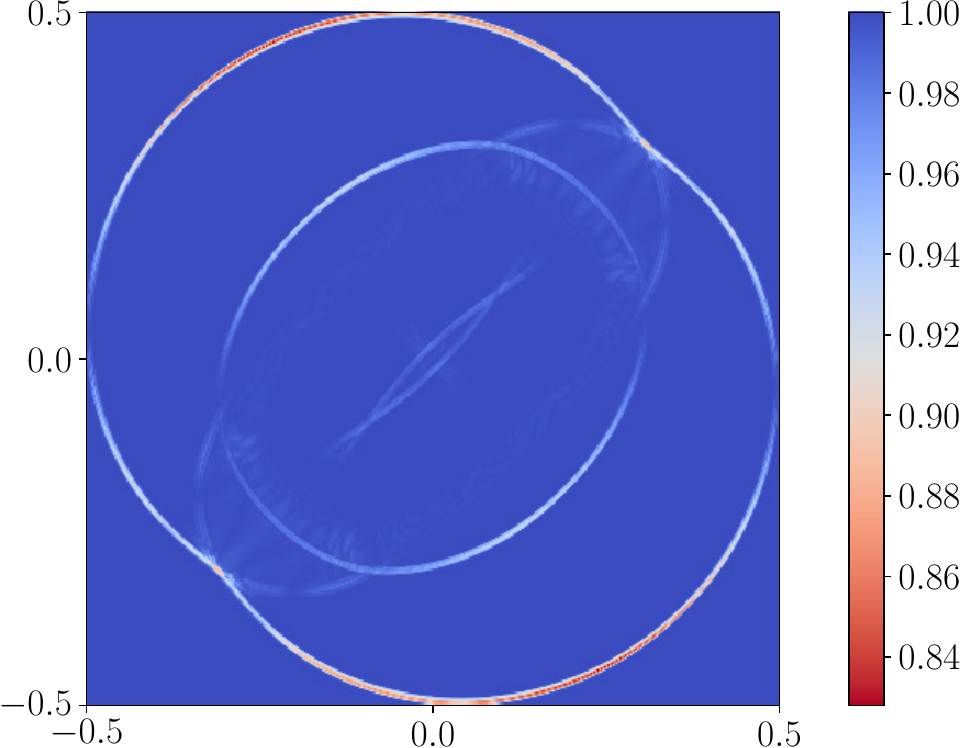}
			\caption{$\theta_{i,j+\frac12}^s$}
		\end{subfigure}
		\caption{Example \ref{ex:2d_blast}.
			$30$ equally spaced contour lines of the numerical solutions obtained by our PP AF scheme, and the blending coefficients in the shock sensor with $\kappa=1$.}
		\label{fig:2d_blast}
	\end{figure}
\end{example}

\begin{example}[Shock-cloud interaction]\label{ex:2d_shock_cloud}
	It is about a strong shock wave interacting with a dense cloud \cite{Dai_1998_simple_JCP,Toth_2000_$B0$_JCP}.
	A planar shock wave moves from $x = 0.6$ to the right, with the left and right states
	\begin{equation*}
		(\rho, \bv, \bB, p) = \begin{cases}
			(3.86859, ~0, ~0, ~0, ~167.345, ~0, ~2.1826182, ~-2.1826182), &\text{if}~x < 0.6, \\
			(1, ~-11.2536, ~0, ~0, ~1, ~0, ~0.56418958, ~0.56418958), &\text{if}~x>0.6.
		\end{cases}
	\end{equation*}
	There is a circular cloud with $\rho=10$ centered at $(0.8, 0.5)$ with a radius of $0.15$.
	The adiabatic index is $\gamma=5/3$ and the final time is $T=0.06$.
	
	The numerical solutions obtained by our PP AF scheme with $400\times400$ cells and $\kappa=1$ are shown in Figure \ref{fig:2d_shock_cloud}.
	The complex structures due to the interaction are captured well without obvious oscillations and they match those in \cite{Dai_1998_simple_JCP,Toth_2000_$B0$_JCP,Wu_2018_provably_SJSC,Liu_2025_Structure_JCP}.
	Note that our limiting based on the shock sensor is locally activated.
	The PP limitings are important in running this test.
    The simulation stops due to negative pressure if they are not used.
	
	\begin{figure}[htb!]
		\begin{subfigure}[t]{0.325\textwidth}
			\centering
			\includegraphics[height=0.8\linewidth]{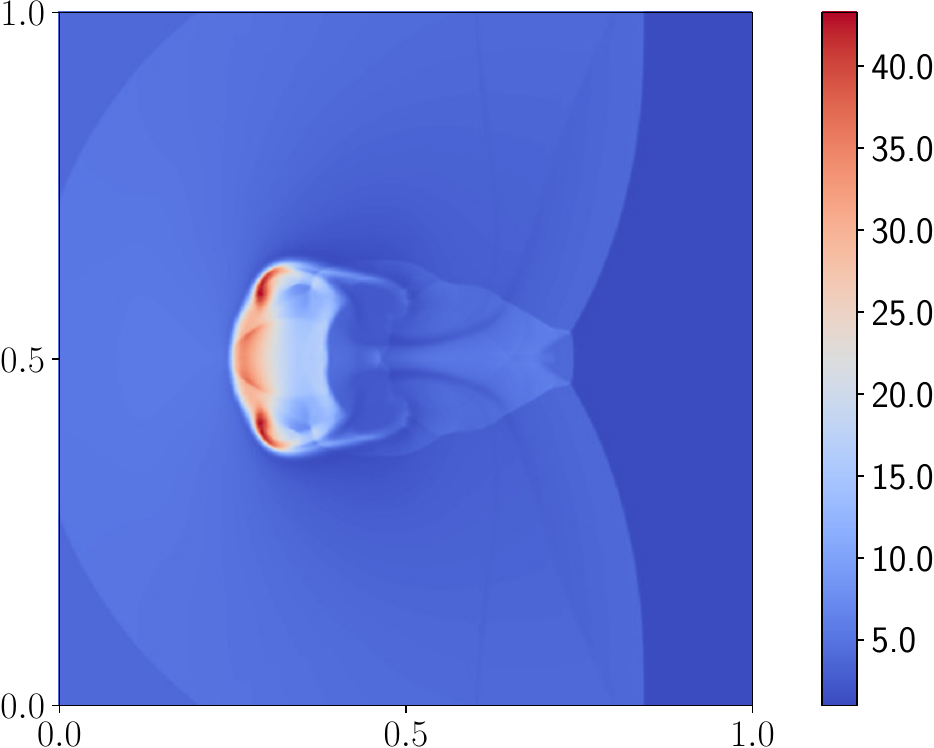}
			\caption{$\rho$}
		\end{subfigure}\hfill
		\begin{subfigure}[t]{0.325\textwidth}
			\centering
			\includegraphics[height=0.8\linewidth]{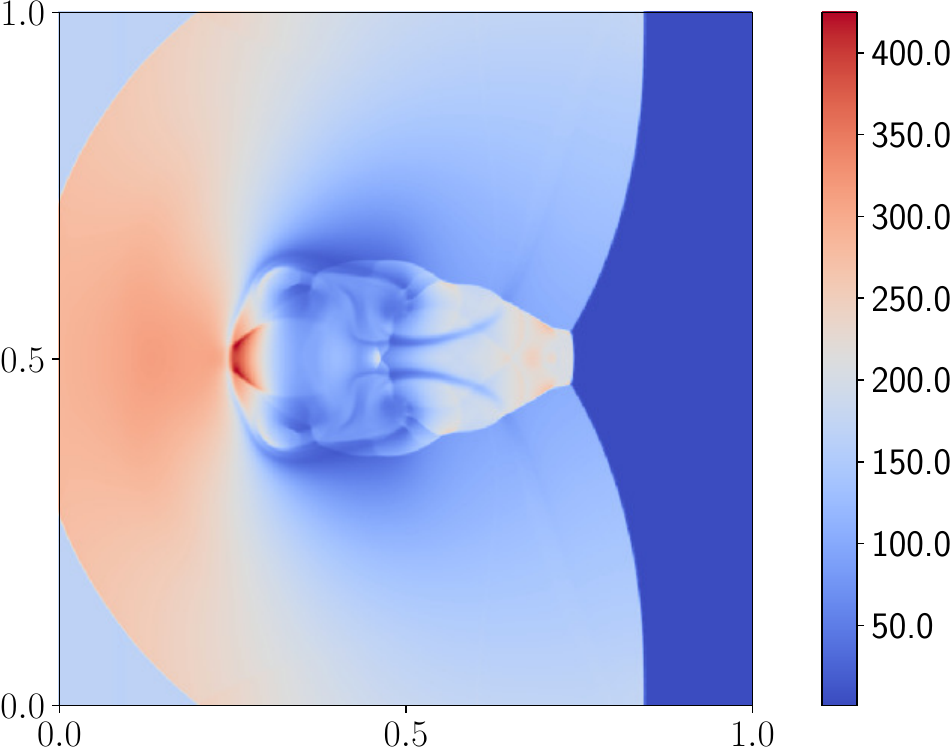}
			\caption{$p$}
		\end{subfigure}\hfill
		\begin{subfigure}[t]{0.325\textwidth}
			\centering
			\includegraphics[height=0.8\linewidth]{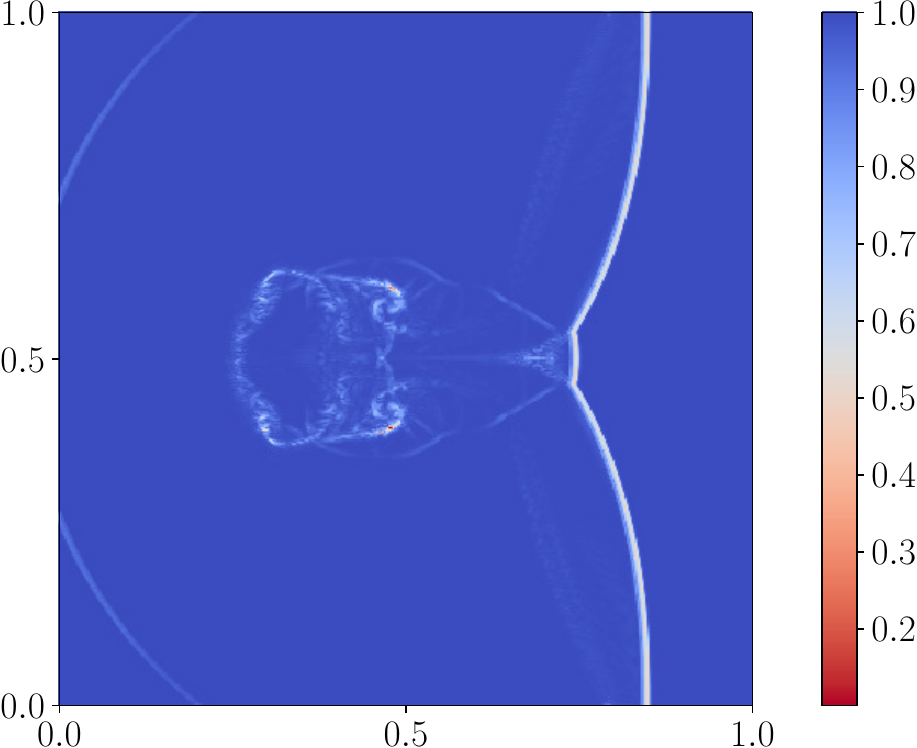}
			\caption{$\theta_{\xr,j}^s$}
		\end{subfigure}
		
		\begin{subfigure}[t]{0.325\textwidth}
			\centering
			\includegraphics[height=0.8\linewidth]{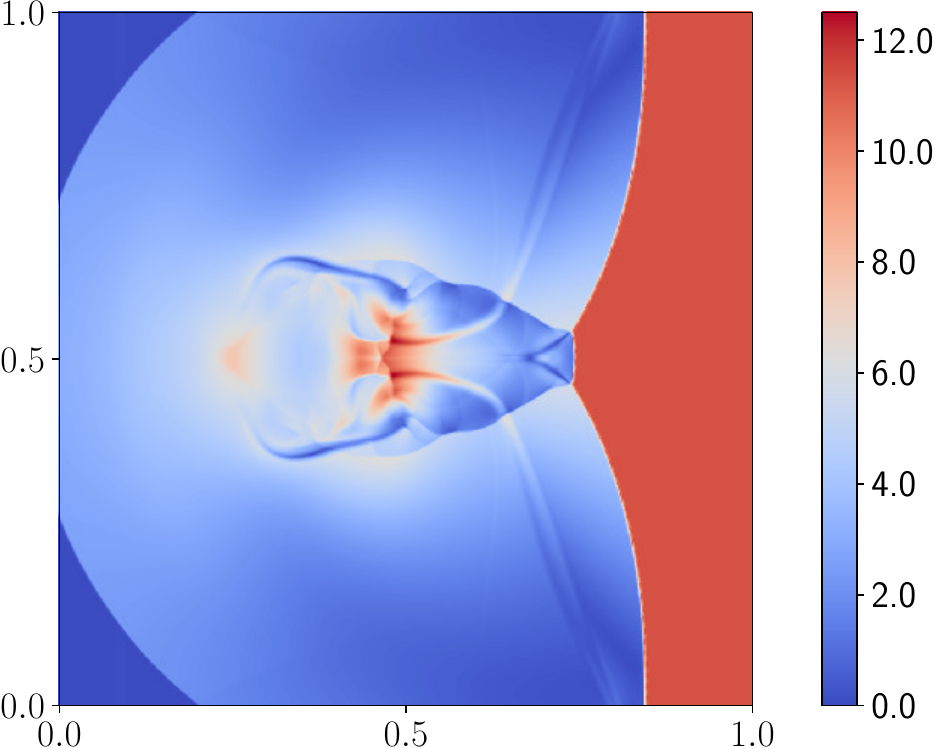}
			\caption{$\norm{\bv}$}
		\end{subfigure}\hfill
		\begin{subfigure}[t]{0.325\textwidth}
			\centering
			\includegraphics[height=0.8\linewidth]{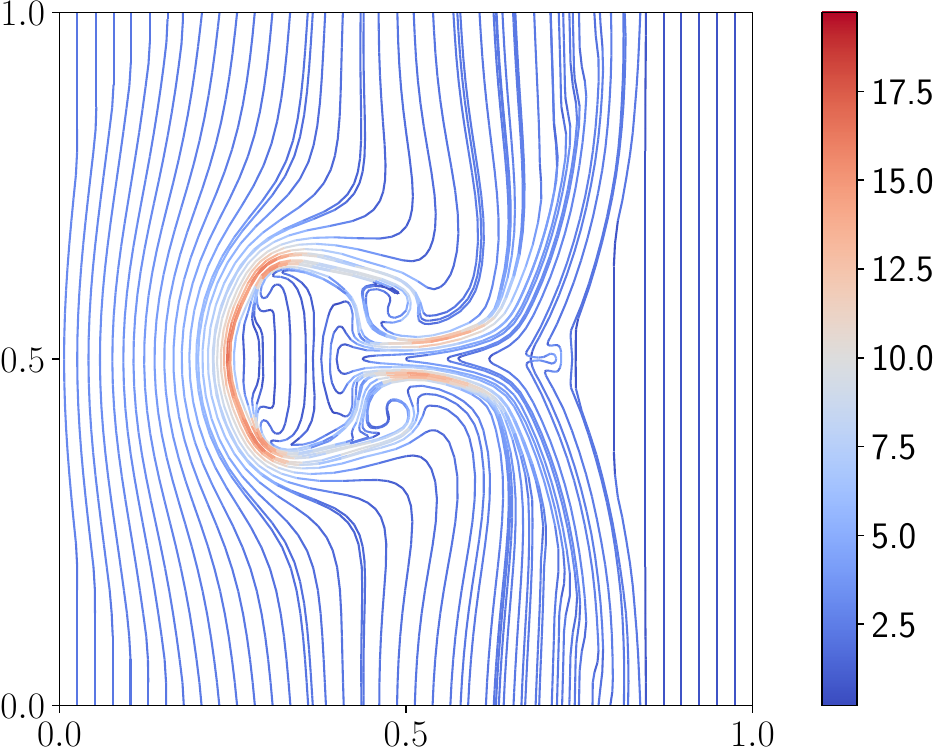}
			\caption{Magnetic lines}
		\end{subfigure}\hfill
		\begin{subfigure}[t]{0.325\textwidth}
			\centering
			\includegraphics[height=0.8\linewidth]{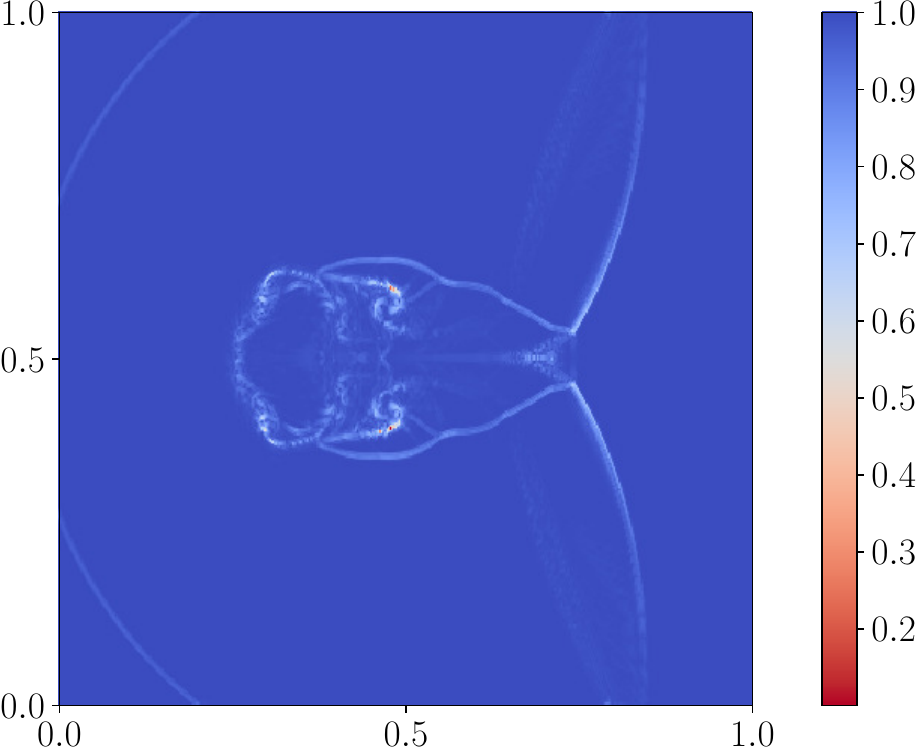}
			\caption{$\theta_{i,\yr}^s$}
		\end{subfigure}
		\caption{Example \ref{ex:2d_shock_cloud}.
			The numerical solutions obtained by our PP AF scheme and blending coefficients in the shock sensor-based limiting with $\kappa=1$.}
		\label{fig:2d_shock_cloud}
	\end{figure}
\end{example}

\begin{example}[MHD jets]\label{ex:2d_jets}
	This is a test problem involving high Mach number jets in a highly magnetized medium \cite{Wu_2018_provably_SJSC}, by adding a magnetic field in the gas dynamical jet in \cite{Balsara_2012_Self_JCP}.
	Following the setup in \cite{Wu_2018_provably_SJSC}, the computational domain is $[-0.5,0.5]\times[0,1.5]$ and the adiabatic index is $\gamma=1.4$.
	The initial ambient fluid is static with $\rho=0.1\gamma$, $p=1$,
	and the magnetic field is initialized in the $y$-direction $B_2=B_a$.
	A jet is injected into the domain with $(\rho, v_1, v_2, p) = (\gamma, 0, 800, 1)$ at the lower boundary when $\abs{x}<0.05$.
	The outflow boundary conditions are applied on other boundaries.
	Here, only the left half is simulated by using a reflective boundary condition at $x=0$.
	The final time is $T=0.002$.
	
	This test problem cannot be simulated without the PP limitings as the solution contains strong discontinuities, low density, and pressure.
	The logarithm of density $\log_{10}\rho$ and pressure $\log_{10}p$ obtained by our PP AF scheme on a $200\times 600$ mesh are shown in Figure \ref{fig:2d_jets}.
	The parameter in the shock sensor is $\kappa=2$ for the three cases $B_a=\sqrt{200}, \sqrt{2000}, \sqrt{20000}$.
	The main flow structures and small-scale features are captured well, comparable to those in \cite{Wu_2018_provably_SJSC}.
	
	\begin{figure}[htb!]
		\begin{subfigure}[b]{0.325\textwidth}
			\centering
			\includegraphics[height=1.0\linewidth]{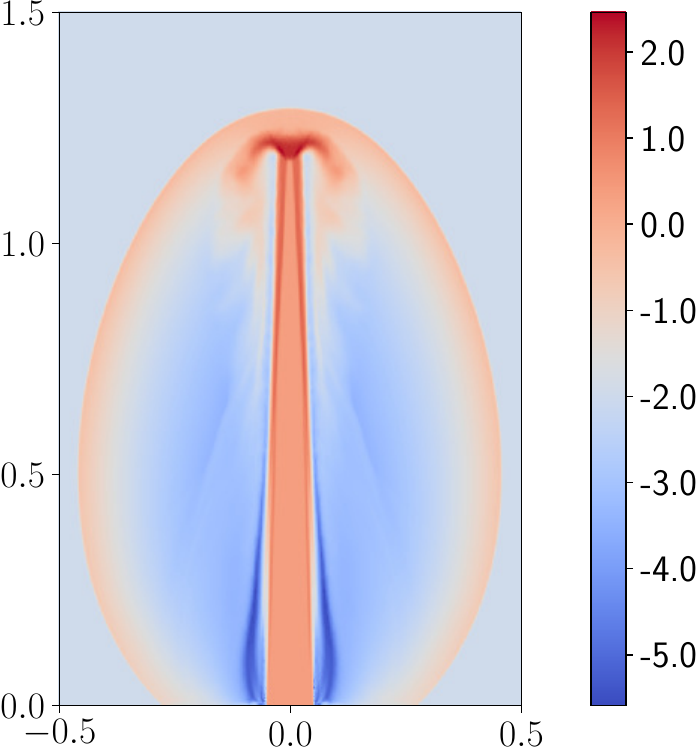}
		\end{subfigure}\hfill
		\begin{subfigure}[b]{0.325\textwidth}
			\centering
			\includegraphics[height=1.0\linewidth]{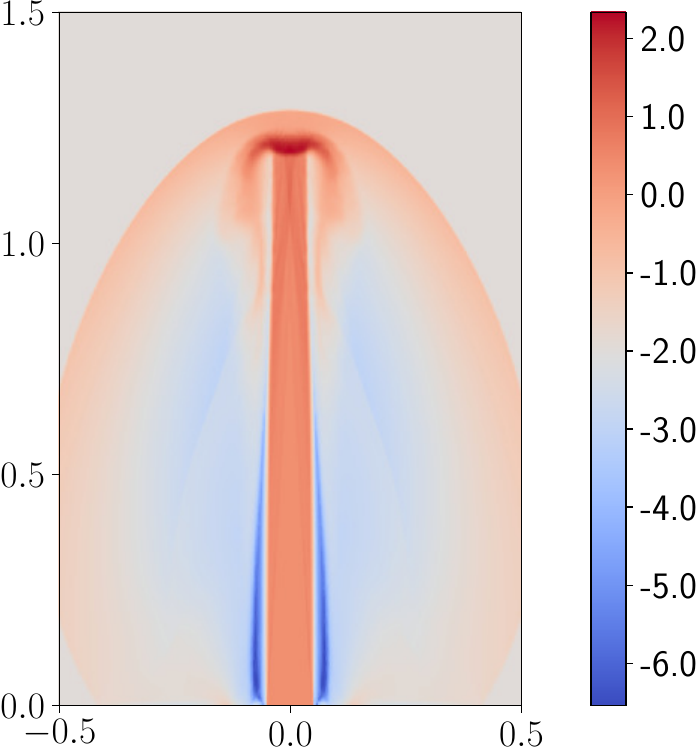}
		\end{subfigure}\hfill
		\begin{subfigure}[b]{0.325\textwidth}
			\centering
			\includegraphics[height=1.0\linewidth]{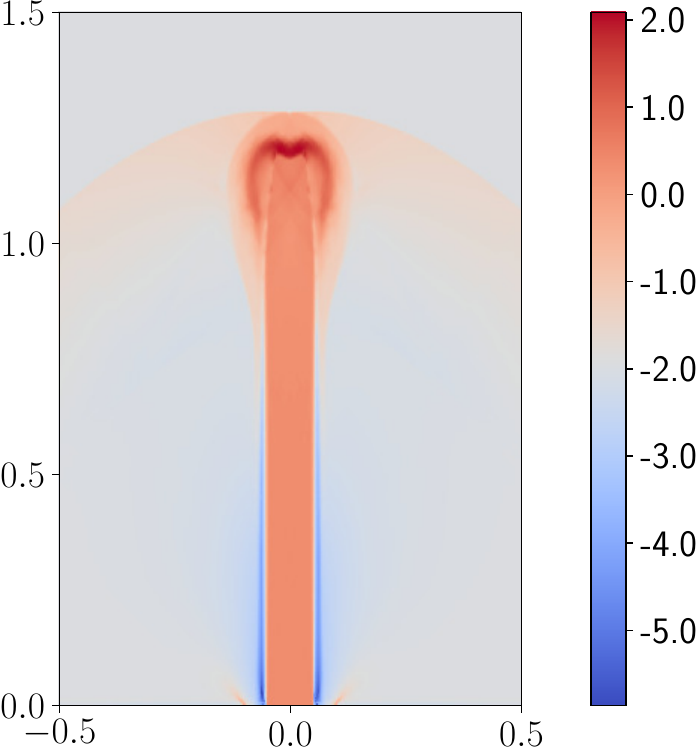}
		\end{subfigure}
		
		\begin{subfigure}[b]{0.325\textwidth}
			\centering
			\includegraphics[height=1.0\linewidth]{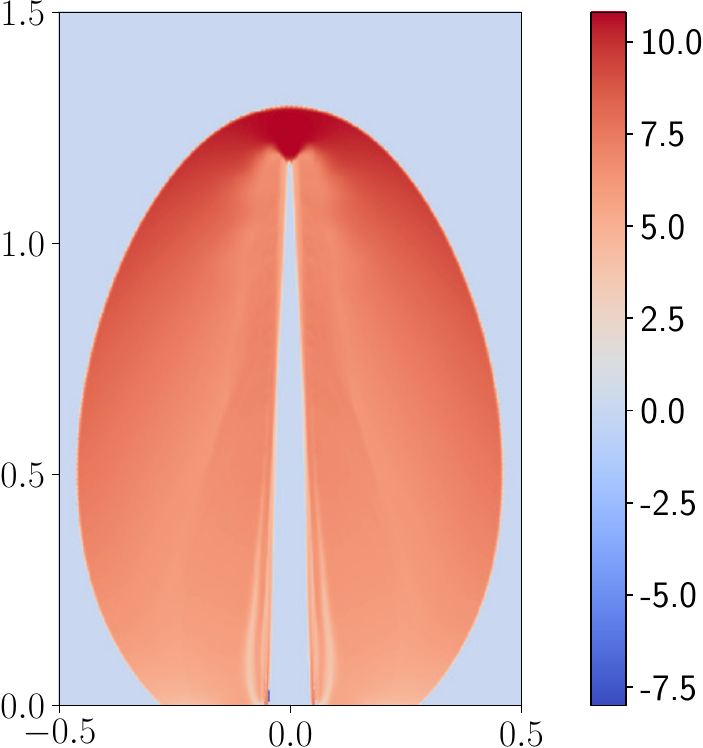}
		\end{subfigure}\hfill
		\begin{subfigure}[b]{0.325\textwidth}
			\centering
			\includegraphics[height=1.0\linewidth]{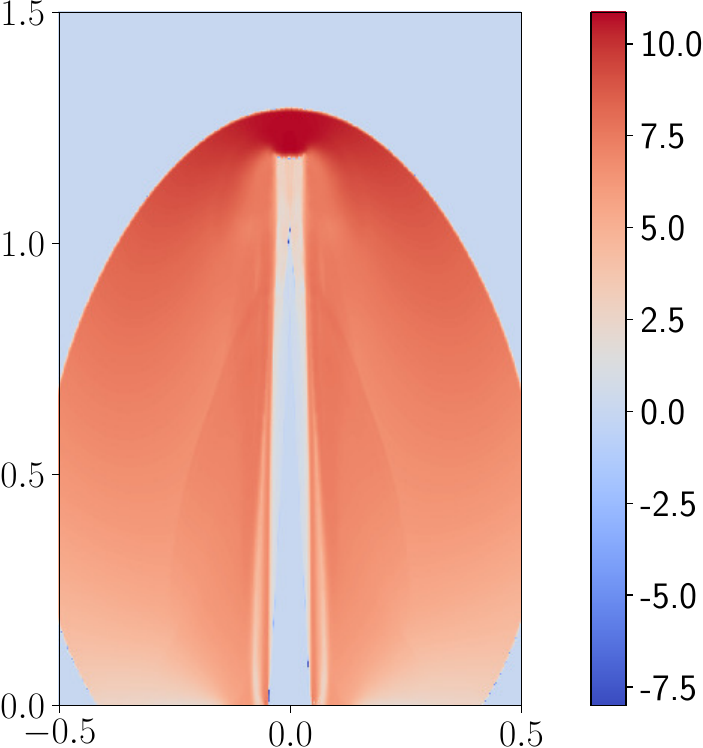}
		\end{subfigure}\hfill
		\begin{subfigure}[b]{0.325\textwidth}
			\centering
			\includegraphics[height=1.0\linewidth]{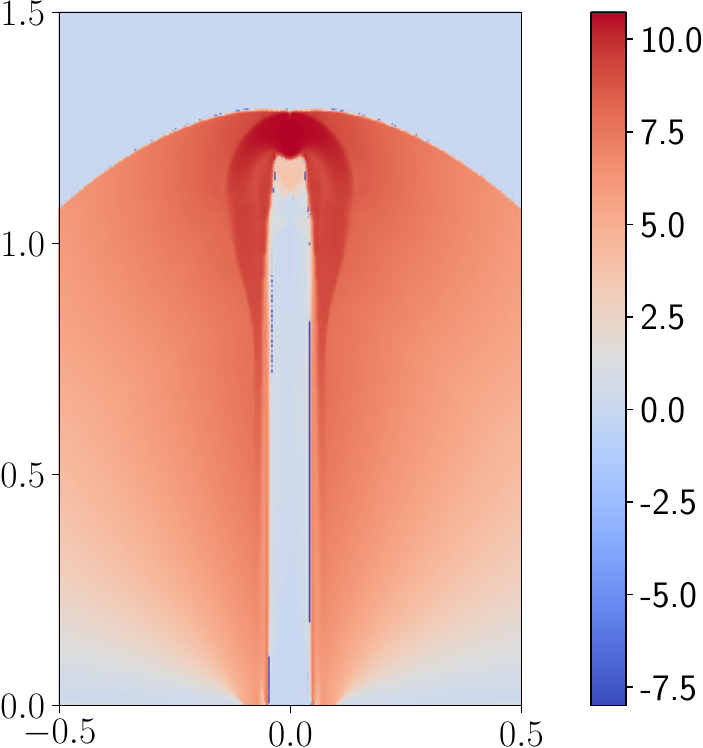}
		\end{subfigure}
		\caption{Example \ref{ex:2d_jets} with $B_a=\sqrt{200}, \sqrt{2000}, \sqrt{20000}$ (from left to right).
			The logarithm of density (top) and pressure (bottom) obtained by our PP AF scheme with $\kappa=2$.}
		\label{fig:2d_jets}
	\end{figure}
	
\end{example}

\begin{remark}
	Due to round-off errors, the limited state may not be PP after the limitings when the scales of the variables differ a lot.
	In this case, the blending coefficients are gradually shrunk by $2^{m}10^{-8}$ with $m=0,1,\dots,9$ until the limited state is PP.
	This case happens rarely, e.g., only $8$ times during $6607$ time steps in the third jet problem with $B_a=\sqrt{20000}$ and $401\times1201$ DoFs.
\end{remark}

\section{Conclusion}\label{sec:conclusion}
This paper has developed a third-order PP AF scheme for solving the ideal MHD equations, with the help of the Godunov-Powell source term to deal with the divergence-free constraint.
The cell average was evolved following the standard finite volume method with a suitable discretization for the nonconservative source term.
This part was free from any Riemann solver due to the continuous representation of the numerical solution at cell interfaces.
The point value update was built on the LLF FVS and a central difference for the source term.
The scheme maintained the compact spatial stencil of the original AF scheme.
The PP limitings for both the cell average and point value were presented to improve robustness for flows containing low density or pressure,
where the parametrized flux limiter and scaling limiter were used to blend the high-order AF scheme and first-order PP LLF scheme, respectively.
To further suppress oscillations, a new shock sensor was employed in the flux limiting. 
Several numerical tests verified the third-order accuracy, PP property, and shock-capturing ability of our scheme.
It was also shown that the Godunov-Powell source term and its suitable discretization played an important role in the control of divergence error and improved stability.

% #############################################################################################################################

\section*{Acknowledgement}

% #############################################################################################################################
JD was supported by an Alexander von Humboldt Foundation Research Fellowship CHN-1234352-HFST-P.
The work of Praveen Chandrashekar was supported by the Department of Atomic Energy, Government of India, under project No.~12-R\&D-TFR-5.01-0520.
CK acknowledged funding by the Deutsche Forschungsgemeinschaft (DFG, German Research Foundation) within \textit{SPP 2410 Hyperbolic Balance Laws in Fluid Mechanics: Complexity, Scales, Randomness (CoScaRa)}, project No.~525941602.
% #############################################################################################################################

% \bibliographystyle{siamplain}
% \bibliography{~/Research/references.bib}

\newcommand{\etalchar}[1]{$^{#1}$}
%\input{PPAF_MHD.bbl}
%\begin{thebibliography}{LMSWZ04}

%\end{thebibliography}

%\appendix

%\section{Section 1} 

% #############################################################################################################################
%Insert your appendix here.
% #############################################################################################################################

\end{document}